\documentclass[letterpaper, 10 pt, conference]{ieeeconf}  % Comment this line out
\usepackage{epsfig}
\usepackage{url}
\usepackage{psfrag}
\usepackage{amsfonts}
\usepackage{graphicx}
\usepackage{framed}
\usepackage{epstopdf}
\usepackage{tikz}
\usepackage{array}
\usepackage{mathtools}
\usepackage{subfig}
\usepackage{steinmetz}
\usepackage{balance}
\usepackage{etex}
\usepackage{url}
\usepackage{color}
\usepackage[usenames,dvipsnames]{pstricks}
\usepackage{pst-grad} % For gradients
\usepackage{pst-plot} % For axes
\IEEEoverridecommandlockouts                              % This command is only
\title{\LARGE \bf
Stability, convergence and bifurcation in some models of chemical kinetics
}

\author{Abuthahir Abdulrahuman, Kalyan Sundar Chakrabarti, and Gaurav Raina% <-this % stops a space
\thanks{Abuthahir Abdulrahuman and Kalyan Chakrabarti are with the School of Interwoven Arts and Sciences, Krea University, Sri City, Andhra Pradesh 517646, India. Gaurav Raina is with the Department of Electrical Engineering, Indian Institute of Technology Madras, Chennai-600036, Tamilnadu,  India. 
        { \tt\small \newline email ids: abuthahir.abdulrahuman@krea.edu.in, kalyan.chakrabarti@krea.edu.in, gaurav@ee.iitm.ac.in}}%
}

\begin{document}

\maketitle
\thispagestyle{empty}
\pagestyle{empty}

%%%%%%%%%%%%%%%%%%%%%%%%%%%%%%%%%%%%%%%%%%%%%%%%%%%%%%%%%%%%%%%%%%%%%%%%%%%%%%%%
\begin{abstract}
In this paper, we analyze the stability, convergence, and bifurcation properties of the Boissonade-De Kepper (BD) model which played a key role in the development of nonlinear chemical dynamics. We first outline conditions for local stability, which may help guide design considerations. Then, we show that the BD model undergoes a Hopf bifurcation when the stability condition gets violated. Using Poincar\'{e} normal forms and center manifold theory, we derive explicit analytic expressions for determining the type of the Hopf bifurcation and the stability of the limit cycles. This provides insights on the system dynamics just beyond the stable regime. Some of the analytical insights are corroborated with numerical computations. We also show that the mathematical results obtained in this paper may have wider applicability beyond the BD model. 
\end{abstract}
%%%%%%%%%%%%%%%%%%%%%%%%%%%%%%%%%%%%%%%%%%%%%%%%%%%%%%%%%%%%%%%%%%%%%%%%%%%%%%%%
\section{INTRODUCTION}
The existence of chaotic dynamics has been noted in the chemical and biochemical systems by multiple investigators ~\cite{Toker2020}. The dynamics of the color change of the Briggs-Rauscher oscillating reaction systems ~\cite{Briggs1973} has been modeled using the Boissonade-De Kepper model ~\cite{boissonade1980}, which serves as one of the illustrative examples of chaotic dynamics in chemical systems. Similar oscillation and bistability is also seen in \emph{vitro} biochemical systems involving the oxidation of the NADH by $O_2$, catalyzed by the horse-radish peroxidase, in a stirred tank reactor where $O_2$ can enter by diffusion from the gas phase ~\cite{Olsen1977}. While these examples use a continuously stirred tank reactor (CSTR) to homogenize the systems, there have been efforts to introduce delayed feedback control in the CSTR ~\cite{Lekebusch1995}.

The use of engineering principles to design biological circuits have been a recent development ~\cite{Niederholtmeyer2015}. The design of the biological systems requires a predictable monotonic or periodic behavior and in order to achieve stability, it is important to understand the nonlinear dynamics and the convergence to predictable behavior of the designed biological systems. This problem is more complex compared to the chemical or biochemical reactions in the CSTR in the absence of the instant or controlled feedback delay. We have analyzed the general class of equations known to exhibit chaotic dynamics for their convergence to stable behavior.

Delay differential equations (DDEs) present both significant opportunities and unique difficulties for chemical modeling. 
There are many reasons why one might wish to expand the class of models available to chemists to include delayed variable formulations. The emphasis in creating a DDE model of a chemical system is shifted from cataloging intermediates and their reactions to describing the dynamic relationships between the concentrations of key species. In many chemical oscillators, it is possible to construct relatively simple models, involving only a handful of key species, that accurately mimic the most important features of the dynamics. As a result, fewer concentration variables will generally appear than in a classical mass-action mechanism.  One aspect that has been remarked upon in a number of mechanisms and models for oscillating chemical reactions is the presence of time-delayed feedback.

The presence of feedback delays makes the system infinite-dimensional and may pose numerous theoretical and practical challenges. In general, the stability of a closed-loop system is sensitive to feedback delays, which normally necessitates a detailed stability analysis. For example, see \cite{cracium2020,insperger2015,kamath2018,raina2005} for some stability and bifurcation analysis of dynamical systems with feedback delays.  Delay dynamical systems are often modeled using delay differential equations to facilitate a mathematical analysis of their performance and dynamics. The initial, and in fact very common, style of stability analysis for nonlinear time-delayed systems is to first linearize the system about its equilibrium and then study the stability properties of the linearized system. However, the feedback delays of a nonlinear dynamical system may result in various complex dynamics like bifurcation, chaos, etc. So, it looks appealing to have an analytical methodology that may allow us to investigate the effect of some nonlinear terms on the system dynamics. Local bifurcation theory is one such methodology \cite{hassard1981}. Moreover, without an understanding of the dynamics of the system in the unstable regime, choosing an operating point close to the boundary of the stable region could be risky. 
A comprehensive understanding of local bifurcation phenomena may help yield insights into the behavior of the system in the unstable regime. Apart from ensuring stability, it is also important to make sure that the system converges quickly to a stable equilibrium.
% For example, see \cite{dubeycnsns2019,huangcnsns2018,khoshcsf2019,novitzkysiads2019} for some stability and bifurcation analysis of dynamical systems with feedback delays. 

In this paper, we conduct the following: local stability, rate of convergence, non-oscillatory convergence, and Hopf bifurcation analyses for the delayed variant of the Boissonade-De Kepper (BD) model \cite{epstein1991}.
% In this paper, we analyze stability, convergence, and Hopf bifurcation properties of the  
% Boissonade-De Kepper (BD) model.
 Our contributions can be summarized as follows.
\begin{enumerate}
\item In the stability analysis, we establish a necessary and sufficient condition to ensure the stable operation of the system. We show that, if the stability condition gets violated, then the system would undergo a Hopf bifurcation, which leads to the emergence of limit cycles. We also derive a sufficient condition for stability. The stability conditions enable us to understand the trade-offs between various system parameters.
\item  We conduct a rate of convergence analysis that enables us to understand the impact of time delay on the convergence rate. We also derive a necessary and sufficient condition that guarantees non-oscillatory convergence to the equilibrium.
\item Using the theoretical frameworks of Poincar\'{e} normal form and the center manifold theorem \cite{hassard1981}, we conduct a detailed Hopf bifurcation analysis that enables us to determine the direction and stability of the emerging limit cycles.
\item To develop a better understanding of how the system dynamics would vary with non-linearity, we also analyze the bifurcation properties of a quadratic model, where the cubic term in the original BD model is replaced by a quadratic term.
\item We show that the BD model, which incorporates cubic control law, can undergo both super-critical and sub-critical Hopf, depending on the parameter values. Whereas, in the case of quadratic model, the Hopf bifurcation is always sub-critical. In general, the occurrence of a sub-critical Hopf is undesirable as it may give rise to either limit cycles with a large amplitude or unstable limit cycles \cite{strogatz2018}. Therefore, our results tend to favor the cubic model. 
\item We also validate some of our analytical insights using numerical simulations and bifurcation diagrams.
\item In the Appendix, we derive a simple closed-form analytic expression for the quantities required to determine the type of the Hopf bifurcation and the stability of the bifurcating periodic solutions of a general first-order non-linear delay differential equation. It is important to highlight that this result has wider applicability beyond the BD model, and can be extended to other non-linear delayed systems as well. To highlight the implications of our general results on Hopf bifurcation, we also apply those results to analyze the bifurcation properties of the Nicholson's Blowflies equation \cite{nicholsonbfref}, which has been extensively used in the context of population dynamics.
\end{enumerate}
The rest of this paper is structured as follows. In Section 2, we outline the models under study. In Section 3, we investigate the local asymptotic stability. The convergence and local Hopf bifurcation analyses are outlined in Sections 4 and 5. Finally, in Section 6, we summarize our key insights and suggest some avenues for further research. For ease of exposition, the Hopf bifurcation analysis is contained in an Appendix.
% 3.1, the bifurcation takes place when  crosses 0 to the right ( > 0), and the bifurcating
% periodic solution is asymptotically stable. We have carried out numerical simulations on
% system (1) using Mathematica with these parameter values, and for different choices of
% initial conditions and for different delays. The simulations support our theoretical result.
%%%%%%%%%%%%%%%%%%%%%%%%%%%%%%%%%%%%%%%%%%%%%%%%%%%%%%%%%%%%%%%%%%%%%%%%%%%%%%%%
\section{MODELS}
\subsection{Motivation}
Complex chemical reactions held far from equilibrium exhibit a variety of phenomena that include multiple stationary states, periodic oscillations, and chaotic oscillations. Nonlinearities in reaction mechanisms may lead to such behavior: these arise from autocatalysis, cross catalysis, and other types of feedback loops. Many chemical processes displaying these complex phenomena are modeled by ordinary differential equations. However, chemical processes with time delays in their mechanism display similar characteristics and are described by delay differential equations (DDEs) in which the rates of change of the variables depend upon their past values. Also, the level of description afforded by a DDE model is often closer to our state of knowledge than is a detailed mechanism in which a certain amount of speculation about intermediate species is a necessary element.

In this paper, we analyze the local stability, convergence, and Hopf bifurcation in a delayed variant of 
Boissonade-De Kepper (BD) model \cite{epstein1991}. Non-linear systems often exhibit periodic oscillations when they lose stability. Researchers have established that the oscillatory behavior observed in many biological, chemical, and engineered systems may be explained by the occurrence of the Hopf bifurcation (Marsden and McCracken, 2012). Local stability and convergence analyses rely mainly on the linear terms. Whereas, in the bifurcation-theoretic analysis, we have to take non-linear terms into consideration, which helps to determine the type of the Hopf bifurcation and the stability of the bifurcating limit cycles. It would be interesting to examine how the nature of the Hopf bifurcation changes with the non-linearity. For example, what happens if the cubic term in \eqref{eq:bd_delayedmodel} is replaced by a quadratic term. At least one clearly motivated design objective would be to choose control laws which not only ensure stability, but also offer better bifurcation-theoretic properties. To that end, we also analyze the bifurcation properties of a quadratic model, which results from replacing the cubic term in the original BD model by a quadratic term.
\subsection{Model description}
The Boissonade-De Kepper (BD) model \cite{boissonade1980} contains a primary variable $x$, whose dynamics are governed by a cubic rate law, and a feedback
variable $y$, which, with appropriate values of the parameters, provides a delayed feedback that causes the primary bistable system to become oscillatory. The model is given by \cite{boissonade1980}
\begin{align}
\dot{x}(t) &= - (x^3 -\mu x + \Lambda) - k y, \label{eq:bd1}\\
\dot{y}(t) &= (x - y)/T \label{eq:bd2}.
\end{align}
The role of the variable $y$ in the above equations is to generate a delayed feedback. Also, at steady state, we have $x = y$. By dropping \eqref{eq:bd2}, and replacing $y(t)$ in \eqref{eq:bd1} by $x(t - \tau)$, we get the delayed variant of the BD model as \cite{epstein1991}
\begin{equation}
 \dot{x}(t) = -(x^3(t) - \mu x(t) + \Lambda) - k x(t-\tau). \label{eq:bd_delayedmodel}
\end{equation}
For our analysis, we consider $k > \mu$ so that the equation \eqref{eq:bd_delayedmodel} has a unique equilibrium. The quadratic version of the above model is
\begin{equation}
 \dot{x}(t) = -(x^2(t) - \mu x(t) + \Lambda) - k x(t-\tau). \label{eq:bd_delayedmodel_quadratic}
\end{equation}
\section{LOCAL STABILITY}
Understanding linear equations can also give us some qualitative insights about a more general non-linear problem. 
% Since we analyze only the local stability, the linearization would give sufficient information to deduce whether or not the system converges to equilibrium. The idea is similar to what we do in calculus in trying to approximate a function by a line with the right slope. 
To linearize the non-linear system, we write the Taylor series expansion of the system about its equilibrium point, and include only the linear terms. In this section, we derive conditions to ensure local asymptotic stability of \eqref{eq:bd_delayedmodel}. We also establish that the system loses local stability via a Hopf bifurcation. We introduce an exogenous non-dimensional bifurcation parameter, $\eta > 0$, to drive the system just into the unstable regime. Let us consider the perturbation $u(t)= x(t) - x_e$, where $x_e$ is the equilibrium which is given by $x_e^3 + (k- \mu) x_e + \Lambda= 0$. Now, the Taylor series expansion of \eqref{eq:bd_delayedmodel}  about the equilibrium ($x_e$) is given by
\begin{equation}
 \frac{d}{dt}u(t) = \eta \big((-3 x^2_e + \mu) u(t) - k u(t - \tau) -3x_e u^2(t) - u^3(t)\big). \label{eq:bd_taylor_exp}
\end{equation}
%where $x_e$ is given by $x_e^3 + (k- \mu) x_e + \Lambda= 0.$
Now, the linearized version of the actual non-linear system is given by 
\begin{equation}
  \frac{d}{dt}u(t) = \eta \big((-3 x^2_e + \mu) u(t)- k u(t-\tau)\big). \label{eq:bd_linearmodel}
\end{equation}
{\color{black}Similarly, the linearized model of \eqref{eq:bd_delayedmodel_quadratic} is given by 
\begin{equation}
  \frac{d}{dt}u(t) = \eta \big((-2 x_e + \mu) u(t)- k u(t-\tau)\big). \label{eq:bd_linearmodel_quadratic}
\end{equation}
From \eqref{eq:bd_linearmodel} and \eqref{eq:bd_linearmodel_quadratic}, we can note that the linearized model of \eqref{eq:bd_delayedmodel_quadratic} is quite similar to that of \eqref{eq:bd_delayedmodel}. Therefore, the results of linear analyses like stability, rate of convergence, and non-oscillatory convergence of \eqref{eq:bd_delayedmodel} can be extended to \eqref{eq:bd_delayedmodel_quadratic}. }

The linearized stability of \eqref{eq:bd_delayedmodel} is given by the stability of the trivial fixed point of \eqref{eq:bd_linearmodel}. The stability of \eqref{eq:bd_linearmodel} is given by the roots of the associated characteristic equation. Looking for exponential solutions, the characteristic equation of (\ref{eq:bd_linearmodel}) is given by
\begin{equation}
 \lambda + \eta a + \eta b e^{-\lambda \tau} = 0, \label{eq:bd_chareq}
\end{equation}
where $a = (3x_e^2 - \mu)$ and $b = k$. Let us consider the case where $\eta, a, b > 0$ and $b > a$.

For the system to be stable, all the roots of the characteristic equation should lie in the left half of the complex plane. For $\tau = 0$, the characteristic equation has a negative real root, and hence the system is asymptotically stable. However, when $\tau > 0$ the roots may cross the imaginary axis for some values of the system parameters, and hence the stability of the system cannot be guaranteed. Therefore, the condition for the crossover defines the bounds on the system parameters to maintain stability. We are interested in finding a critical value at which a root of this equation transitions from having negative to having positive real parts. If this is to occur, there must be a boundary case, such that the characteristic equation has a purely imaginary root. Therefore, to find the critical condition, we substitute $\lambda = \pm j\omega$, $\omega>0$ in (\ref{eq:bd_chareq}). Then, we break the polynomial up into its real and imaginary parts, and write the exponential in terms of trigonometric functions to obtain
% 
% 
% With $\tau = 0$, the characteristic equation has a negative real root, and hence the system is asymptotically stable. For $\tau > 0$, let $\lambda = j\omega$, $\omega > 0$ which gives
\begin{eqnarray}
 \eta a + \eta b \cos(\omega \tau) &= 0 \label{eq:char_real}\\
 \omega - \eta b \sin(\omega \tau) &= 0. \label{eq:char_img}
\end{eqnarray}
For $\omega > 0$, we get $\cos(\omega \tau) < 0$  and $\sin(\omega \tau > 0)$, giving
\begin{equation}
 2n\pi + \pi/2 < \omega \tau < 2n \pi + \pi, \quad \ n = 0,1,2,....
\end{equation}
We only treat the case $n=0$. Solving \eqref{eq:char_real} and \eqref{eq:char_img}, we obtain
\begin{eqnarray}
 \omega_0 = \eta_c \sqrt{b^2 - a^2}\\
 \eta_c \tau \sqrt{b^2 - a^2} = \cos^{-1}(-a/b)
\end{eqnarray}
where $\eta_c$ denotes the critical value of $\eta$ at $\omega = \omega_c$. To show that the system undergoes a Hopf bifurcation at $\eta_c$, we need to satisfy the following transversality condition of the Hopf spectrum \cite{hassard1981}\\
$$\mathbf{Re}\left(\dfrac{d\lambda}{d\eta}\right)_{\eta=\eta_c}\neq 0.$$\\
In other words, for the occurrence of the Hopf bifurcation, the roots of the characteristic equation should cross the imaginary axis from left to right with non-zero speed. Differentiating equation (\ref{eq:bd_chareq}) with respect to $\eta$, we obtain 
\begin{equation}
 \left. \dfrac{d\lambda}{d\eta}\right\lvert_{\eta=\eta_c} = \left. \dfrac{-(a+b e ^{-\lambda \tau})}{1 - \eta b \tau e^{-\lambda \tau}}\right\lvert_{\eta=\eta_c.}
%  \left. \dfrac{-(e^{-\lambda\tau_1}+e^{-\lambda\tau_2})}{(2-\rho  (\tau_1e^{-\lambda\tau_1}+\tau_2e^{-\lambda\tau_2})\big)}\right\lvert_{\rho=\rho_c.}
\end{equation}
From the above equation, we get
\begin{equation*}
 \mathbf{Re}\left(\dfrac{d\lambda}{d\eta}\right)_{\eta=\eta_c} = \frac{\eta_c \tau (b^2 - a^2)}{1 + 2 \eta_c a \tau + \eta_c^2 b^2 \tau^2} \  > \ 0,
\end{equation*}
% $\mathbf{Re}\left(\dfrac{d\lambda}{d\rho}\right)_{\rho=\rho_c}=\dfrac{\pi \sin(\omega_0\tau_1)}{A^2+B^2} > 0,$
% where 
% \begin{eqnarray}
% % \omega_0 &=& \pi/(\tau_1+\tau_2),\nonumber\\
% A &=& 1-  \dfrac{\omega_0\cos(\omega_0\tau_1)(\tau_1-\tau_2)}{2\sin(\omega_0\tau_1)},\nonumber\\
% B &=& \pi/2.\nonumber
% \end{eqnarray}
Hence, the system undergoes a Hopf bifurcation at $\eta = \eta_c$, with period $2\pi/\omega_0$. Thus, the \textit{necessary and sufficient} condition for local asymptotic stability of \eqref{eq:bd_delayedmodel} is 
\begin{equation}
\label{eq:ns_cond_1}
\eta_c \tau < \frac{\cos^{-1}(-a/b)}{\sqrt{b^2 - a^2}},
% \rho \frac{(\tau_1+\tau_2)}{2}\sin\left(\dfrac{\pi\tau_1}{\tau_1+\tau_2}\right)<\pi/2.
\end{equation}
where $a = (3x_e^2 - \mu)$ and $b = k$. From \eqref{eq:ns_cond_1}, we can deduce that longer delays increase the region of parameter
space in which oscillatory behavior may occur. 
% Substituting the values of $a$ and $b$ in \eqref{eq:ns_cond_1}, we get
% \begin{equation}
%  \label{eq:ns_cond_2}
% \eta_c \tau < \frac{\cos^{-1}(-\frac{3x_e^2 - \mu}{k})}{\sqrt{k^2 - (3x_e^2 - \mu)^2}}
% \end{equation}
%  $a = (3x_e^2 - \mu)$ and $b = k$.
%  Substituting the value of $\myk$ in ~(\ref{eq:stability_relation}), we get
% \begin{equation}
% \label{eq:nsstability_relation}
%  a\left[1+\frac{2R^*}{C}\right]\sin\left(\dfrac{\pi\tau_1}{\tau_1+\tau_2}\right)<\pi/2.
% \end{equation}

%%%%%%%%%%%%%%%%%%%%
% A few comments are in order. First, i
It is to be noted that the system becomes unstable when the very first conjugate pair of characteristic roots cross the imaginary axis. As the derivative $\mathbf{Re}\left({d\lambda}/{d\eta}\right)_{\eta=\eta_c}$ is positive, the system cannot regain its local stability with further increase in the value of the bifurcation parameter ($\eta$). In other words, an increase in the value of $\eta$ results in the characteristic roots or eigenvalues moving to the right in the complex plane, thereby making it impossible to restore lost stability.
% Next, we can note that, apart from $\eta$, the time  delays ($\tau_1, \tau_2$) also play a role in ensuring local stability. This implies that an appropriate variation in the values of $\tau_1$ and $\tau_2$ can also induce instability. Therefore, the feedback delays can also be used as the bifurcation parameter. 
Further, note that $f_0 = \omega_0/2\pi$, represents the frequency of the bifurcating periodic oscillations.\\

\textit{Sufficient condition}. We now use Nyquist stability criterion to derive sufficient condition for local stability. From the characteristic equation \eqref{eq:bd_chareq}, we obtain the loop transfer function as 
% \begin{equation}
%  \lambda \Bigg(1+ \frac{\eta a }{\lambda}+ \frac{\eta b e^{-\lambda\tau_2}}{\lambda}\Bigg)=0.
%  %\label{eq:apndxce2}
%  \end{equation}
% Now, the loop transfer function is given 
\begin{equation}
L(\lambda)= \frac{\eta b e^{-\lambda\tau}}{\lambda + \eta a}. 
\label{eq:apndxlooptf}
\end{equation}
The next step is to obtain the crossover frequency at which $\phase{L(j\omega)} =\pi$. At this frequency, the magnitude of the loop transfer function should be less than 1, i.e., $\left|L(j\omega)\right| < 1$. Now, substituting $\lambda=j\omega$ in \eqref{eq:apndxlooptf} yields
\begin{equation}
 L(j\omega) = \frac{\eta b e^{-j\omega \tau}}{j\omega + \eta a}.
\end{equation}
Equating $\phase{L(j\omega)}$ to $\pi$, we get
\begin{equation}
 \tan(\omega \tau) = \frac{-\omega}{\eta a}.
 \label{eq:apndx_angle_condn}
\end{equation}
Similarly, the magnitude condition $\left|L(j\omega)\right| < 1$ can be written as 
\begin{equation}
{
 |L(j\omega)| = \frac{\eta b}{\sqrt{\eta^{2} a^{2} + \omega^{2}}} < 1.
 \label{eq:apndx_mag_condn1}
}
\end{equation}
Substituting \eqref{eq:apndx_angle_condn} in \eqref{eq:apndx_mag_condn1}, we obtain
% \begin{equation}
%  \frac{\eta a\sin(\omega_c\tau_1)}{\omega_c} + \frac{\eta b\sin(\omega_c\tau_2)}{\omega_c} < 1.
%  \label{eq:apndx_mag_condn2}
% \end{equation}
% We can rewrite \eqref{eq:apndx_mag_condn2} as
\begin{equation}
 \eta b \tau\frac{\sin(\omega\tau)}{\omega\tau} < 1.
 \label{eq:apndx_mag_condn3}
\end{equation}
From \eqref{eq:apndx_angle_condn} and \eqref{eq:apndx_mag_condn3}, we gather that $\tan(\omega \tau) < 0$ and $\sin(\omega \tau) > 0$. This implies that $\pi/2 < \omega \tau < \pi$. Therefore, the function on the left-hand side of \eqref{eq:apndx_mag_condn3} attain its maxima at $\omega \tau = \pi/2$, and consequently $\sin(\omega \tau) = 1$. This yields the sufficient condition for local stability as 
\begin{equation}
 \eta b \tau < \pi/2
 \label{eq:suff_condn_nyquist}
\end{equation}
where $b = k$. In fact, for $a, b > 0$ and $b>a$, one can show that the minimum value of right-hand side of \eqref{eq:ns_cond_1} is $\pi/2$. Observe that the sufficient condition \eqref{eq:suff_condn_nyquist} does not involve parameters $\mu$ and $\Lambda$. Therefore, if we ensure $\eta b \tau < \pi/2$ then \eqref{eq:bd_delayedmodel} would be locally asymptotically stable, regardless of the other parameter values. 
%the necessary and sufficient  condition for stability 
% Equation (\ref{eq:img0_tcp}) can be written as 
% \begin{equation}
% \eta \eta a_{1}\tau_1\frac{\sin(\omega\tau_1)}{\omega\tau_1}+ \eta \eta a_{2}\tau_2\frac{\sin(\omega\tau_2)}{\omega\tau_2} =1. \label{eq:img1_tcp}
% \end{equation}
% Conditions for local stability yield bounds on the system parameters and feedback delay to ensure locally stable operation. 
From the stability conditions, we can understand that there exist \textit{trade-offs} between the system parameters and the time delay ($\tau$)  for local stability. Therefore, the above insight may be used to guide design of system parameters such that system stability is ensured. 
\section{CONVERGENCE}
In the local stability analysis, we derived stability conditions which enable us to understand the role of various system parameters in ensuring local stability. Now, within the stable regime, it is also important to study the impact of systems parameters on the convergence characteristics of the system. 

\subsection{Rate of convergence}
Rate of convergence is an important performance metric that dictates the time a dynamical system takes to equilibrate, when perturbed. 
%In the context of a transportation network, it is related to the time required to attain
In this subsection, following the style of analysis outlined in \cite{brauer79}, we conduct a rate of convergence analysis for \eqref{eq:bd_delayedmodel}. To do so, it is sufficient to solve the characteristic equation of \eqref{eq:bd_delayedmodel} whose roots determine the convergence characteristics of solutions of \eqref{eq:bd_delayedmodel} completely.
%and study the dependence of convergence characteristics on time delay $\tau$ and parameter $b$ i.e., $k$. 
The analytical results enable us to investigate the impact of various system parameters on the rate of convergence to the equilibrium. Here, we consider $\eta=1$, to get back the original system. Now, the linearized version of the actual non-linear system is given by 
\begin{equation}
  \frac{d}{dt}u(t) = -a u(t)- b u(t-\tau), \label{eq:bd_linearmodel_roc}
\end{equation}
where $a = (3x_e^2 - \mu)$ and $b = k$.

% In this sub-section, we study the dependence of convergence characteristics on time delay $\tau$ and parameter $b$ i.e., $k$. To do so, it is sufficient to solve the characteristic equation of \eqref{eq:bd_delayedmodel} whose roots determine the convergence characteristics of solutions of \eqref{eq:bd_delayedmodel} completely. We closely follow the style of analysis outlined in \cite{brauer79}.
To analyze the dependence of convergence characteristics on the time delay $\tau$, we require $b \neq 0$. Recall that the characteristic equation is given by
\begin{equation}
 \lambda +  a +  b e^{-\lambda \tau} = 0. \label{eq:bd_chareq_roc}
\end{equation} 
Using $\lambda \tau = z$, $-a \tau = p$ and $-b \tau = q$, the characteristic equation \eqref{eq:bd_chareq_roc} becomes 
\begin{equation}
 (p - z)e^{z} + q = 0. \label{eq:bd_chareq_roc_mod}
\end{equation}
If $-\alpha < 0$, is the real part of a root of \eqref{eq:bd_chareq_roc_mod}, then \eqref{eq:bd_linearmodel_roc} has a solution of the form $e^{(-\alpha/\tau)t}$, which is a decaying function of $t$. If every root of \eqref{eq:bd_chareq_roc_mod} lies in the left half of the complex plane, the solution of \eqref{eq:bd_linearmodel_roc} is asymptotically stable. Here, $\alpha/\tau$ is the convergence rate, the rate at which the stable system approaches equilibrium. The inverse of the rate of convergence, $\tau/\alpha$ is the characteristic return time, where $-\alpha$ is the largest of the real parts of all the roots of \eqref{eq:bd_chareq_roc_mod}. A necessary and sufficient condition for all roots of \eqref{eq:bd_chareq_roc_mod} to be in the left half of the complex plane, given in \cite{hayes50} is,
\begin{equation}
 p < 1, \quad p < -q < \frac{u_1}{\sin(u_1)}, \label{eq:hayes_cond1}
\end{equation}
where $u_1$ is the solution of the equation
\begin{equation}
 u = p \tan(u) \label{eq:u1},
 \end{equation}
in $u \in (0, \pi)$, with $u_1 = \pi/2$ if $p = 0$. We now perform a change of variable $z = \chi - \sigma \tau$ , where $\sigma = \alpha/\tau$ and transform the characteristic equation \eqref{eq:bd_chareq_roc_mod} to 
\begin{equation}
(p + \sigma \tau - \chi) e^{\chi} + q e^{\sigma \tau} = 0. \label{eq:bd_chareq_roc_mod_1}
\end{equation}
Let $\sigma$ be the supremum of the solution of \eqref{eq:bd_chareq_roc_mod_1} over
$(0, \infty)$ which guarantees that all roots of the transformed characteristic equation \eqref{eq:bd_chareq_roc_mod_1} lie in the left half of the complex plane, and then $\sigma$ is the rate of convergence to the equilibrium of \eqref{eq:bd_linearmodel_roc}. For the transformed characteristic equation \eqref{eq:bd_chareq_roc_mod_1}, the set of inequalities \eqref{eq:hayes_cond1} can be restated as
\begin{eqnarray}
 p + \sigma \tau &=& (-a + \sigma)\tau < 1, \label{eq:hayes_cond2}\\
  p + \sigma \tau &=& (-a + \sigma)\tau < -q e^{\sigma \tau} = b\tau e^{\sigma \tau},\label{eq:hayes_cond3}\\
   -q e^{\sigma \tau} &=& b\tau e^{\sigma\tau} < \frac{u_2}{\sin(u_2)},\label{eq:hayes_cond4}
\end{eqnarray}
where $u_2$ is the solution of the equation
\begin{equation} 
u = (-a + \sigma) \tau \tan(u), \label{eq:bd_roc_u2}
\end{equation}
in $u \in (0, \pi)$, with $u_2 = \pi/2$ if $ -a + \sigma = 0$. We choose $\sigma$ to be the supremum of the solution of \eqref{eq:bd_chareq_roc_mod_1} over $(0, \infty)$ satisfying \eqref{eq:hayes_cond2} - \eqref{eq:hayes_cond4} inequalities. Equation \eqref{eq:hayes_cond3}
can be re-written as
\begin{equation}
 (-a + \sigma)\tau e^{(a-\sigma)\tau} < b\tau e^{a\tau}.\label{eq26}
\end{equation}
Consider the function,
\begin{equation}
 g(u) = \frac{u}{\sin(u)} e^{u/\tan(u)}.
\end{equation}
In $u \in (0,\pi)$, $u/\sin(u)$ is an increasing function of $u$, $u/\tan(u)$ is a decreasing function of $u$. So, $g(u)$ is an increasing function of $u$. It is easy to observe that
$g(0) = 1/e,\ g(\pi/2) = \pi/2,\ \lim_{u \rightarrow \pi} g(u) = \infty$. Using \eqref{eq:bd_roc_u2}, the inequality \eqref{eq:hayes_cond4} can be written as
\begin{equation}
 b\tau e^{a\tau} < g(u_2).\label{eq27}
 \end{equation}
It can be observed from \eqref{eq:bd_roc_u2} that $u_2$ is a decreasing function of $\sigma$.  To obtain the maximum $\sigma$ satisfying the inequality \eqref{eq27}, we need to solve its corresponding equality. Similarly, $(-a+\sigma)\tau$ from \eqref{eq:hayes_cond2} and $(-a + \sigma)\tau e^{(a -\sigma)\tau}$ from (\ref{eq26}) are increasing functions of $\sigma$, and to obtain the maximum $\sigma$ satisfying the inequalities (\ref{eq:hayes_cond2}), (\ref{eq26}), we need to solve the corresponding equalities. If there is no solution $\sigma$ satisfying the
equalities corresponding to the inequalities (\ref{eq:hayes_cond2}), (\ref{eq26}) and (\ref{eq27}), there is no restriction on $\sigma$. In summary, these results are as follows:
Let $\sigma_1, \sigma_2, \sigma_3$ be the solutions of
\begin{eqnarray}
 (-a + \sigma)\tau &=& 1,\label{eq28}\\
(-a + \sigma)\tau e^{(a - \sigma)\tau} &=& b\tau e^{a\tau},\nonumber \\
g(u_2) &=& b\tau e^{a\tau},\label{eq29}\\
u_2 &=& (-a + \sigma)\tau \tan(u_2)\label{eq30},
\end{eqnarray}
respectively. $\sigma_1 = \infty$, for $i = 1, 2, 3$ if the corresponding equality has no solution. Then, the rate of convergence of \eqref{eq:bd_linearmodel_roc} is given by
\begin{equation}
 \sigma = \min[\sigma_1, \sigma_2, \sigma_3]. 
\end{equation}

We first characterize the dependence of the rate of convergence on the time delay $\tau$ for $a, b > 0$. For this analysis, we consider the coefficients $a$ and $b$ of \eqref{eq:bd_linearmodel_roc} to be constant. From $p < 1$ of (\ref{eq:hayes_cond1}), we have $-a \tau < 1$, i.e. $\tau > -1/a$, a trivial condition for $a > 0, \tau \geq 0$. From
$p < -q$ of (\ref{eq:hayes_cond1}), we have $-a < b$ which is also a trivial condition for $a, b > 0$. Using \eqref{eq:u1}, $p < u_1/\sin(u1)$ can be written as $\cos(u_1) < 1$, a trivial condition and $-q < u_1/\sin(u_1)$ can be written as $\cos(u_1) < -a/b$, which is a more stricter condition. So, for the case $a, b > 0$, the only relevant stability condition \eqref{eq:hayes_cond1} is $-q < u_1/\sin(u_1)$. As $(-a + \sigma)\tau e^{(a-\sigma)\tau}$ has a maximum $1/e$ at $(-a + \sigma)\tau = 1$ and $g(u)$ has a minimum value of $1/e$ at $u = 0$, no solution exists for equation (\ref{eq29}) if $b \tau e^{a\tau} > 1/e$ and for equation (\ref{eq30}) if $b\tau e^{a\tau} < 1/e$. For $b > 0$, $b\tau e^{a\tau}$ monotonically increases with $\tau$, having a minimum $(-b/a)(1/e) < 1/e$ for $a, b > 0$ by \eqref{eq:hayes_cond1}, attained at $\tau = -1/a$. Thus, there exists $\tau^{*} > -1/a$ such that 
\begin{equation}
 b \tau^{*} e^{a\tau^{*}} = \frac{1}{e}. \label{eq31}
\end{equation}
For $\tau > \tau^{*}$, no solution exists for \eqref{eq29} and for $0 \leq \tau < \tau^{*}$, equation \eqref{eq30} has no solution. For $0 \leq \tau < \tau^{*}$, let $\sigma_2$ be the solution of equation \eqref{eq29}. Differentiation of \eqref{eq29} with respect to $\tau$ results in 
\begin{equation}
 \frac{d\sigma}{d\tau} = \frac{b\sigma e ^{\sigma \tau}}{1 - b \tau e^{\sigma \tau}}. \label{eq32}
\end{equation}
From \eqref{eq29}, $b e^{\sigma \tau} = -a + \sigma$. Using this, \eqref{eq32} can be rewritten as 
\begin{equation}
 \frac{d\sigma}{d\tau} = \frac{b\sigma e ^{\sigma \tau}}{1 - (-a + \sigma) \tau }. %	\label{eq33}
\end{equation}
At $\tau = 0$, $\sigma_1 = \infty$ and $\sigma_2 = a + b$. The derivative $d\sigma_2/d\tau > 0$ if $(-a +\sigma_2)\tau < 1$, i.e. $\sigma_2 < \sigma_1$. At $\tau = \tau^{*}$ 
\begin{equation}
(-a + \sigma_2)\tau e^{(a - \sigma_2)\tau} = b\tau e^{a \tau} = \frac{1}{e}.
\end{equation}
From the above equation, $(-a + \sigma_2)\tau = 1$, which implies $\sigma_1 = \sigma_2$ at $\tau = \tau^{*}$. For $0 \leq \tau < \tau^{*}$, $\sigma_2 < \sigma_1$ and $d\sigma_2/d\tau > 0$. So, for $0 \leq \tau < \tau^{*}$, the rate of convergence $\sigma = \min[\sigma_1, \sigma_2] = \sigma_2$ increases as $\tau$ increases.

For $\tau > \tau^{*}$, from the first condition of (\ref{eq30}), an increase in $\tau$ results in an increase in $u$, thereby resulting in decrease of $\sigma$ as $\tau$ is increased, observable from second condition of (\ref{eq30}). This implies decrease in the rate of convergence for $\tau > \tau^{*}$. Since $u_2 > 0$ for $\tau > \tau^{*}$, $u_2/\tan(u_2) < 1$, and $(-a + \sigma_3)\tau < 1$. Thus, $\sigma_3 < \sigma_1$.

Result: For $a, b > 0$, the rate of convergence is a monotonically increasing function of $\tau$ given by (\ref{eq29}) in the interval $0 \leq \tau < \tau^{*}$, a monotonically decreasing function of $\tau$ given by (\ref{eq30}) in the interval $\tau > \tau^{*}$, $\tau^{*}$ can be obtained from (\ref{eq31}). {\color{black}See Figure 1 for a graphical representation of the dependence of the convergence characteristics on the time delay $\tau$}.

\begin{figure}[hbtp!]
\vspace{-10mm}
\centering
\psfrag{R}{\hspace{-1.8cm} \small Rate of convergence, $\sigma$}
\psfrag{a}{{\hspace{-0.8cm} \small Time delay, $\tau$}}
\psfrag{tau=1}{\hspace{-2mm}\small{$\tau=1$}}
\psfrag{tau=2}{\hspace{-2mm}\small{$\tau=2$}}
\psfrag{unstable}{\small }
\psfrag{roc}{\small Hopf condition}
\psfrag{0}{\footnotesize{$0$}}
\psfrag{0.0000000}{\hspace{0.4cm}\footnotesize{$0$}}
\psfrag{0.387}[][]{\small{$\tau^{*}$}}
\psfrag{0.7853982}{\hspace{0.3cm}\small{$\pi/4$}}
\psfrag{1.5707963}{\hspace{0.3cm}\small{$\pi/2$}}
\psfrag{0.000}{\hspace{0.10cm}\small{$0$}}
\psfrag{0.0}{\hspace{0.0cm}\small{$0$}}
\psfrag{0.9}{\hspace{-0.350cm}\footnotesize{$(a + b)$}}
\psfrag{2.5}[][]{\footnotesize{$(a + 1/\tau^{*})$}}
\psfrag{1.57}{\small{$\pi/2$}}
\includegraphics[width=0.5\textwidth, height= 5cm]{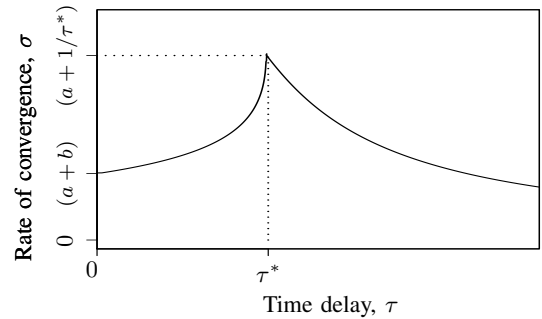}
\caption{Rate of convergence to equilibrium for \eqref{eq:bd_delayedmodel}. Observe that the rate of convergence increases with $\tau$ and reaches maxima of $(a + 1/\tau^{*})$ at $\tau=\tau^{*}$, where $\tau^{*}$ is given by $b \tau^{*} e ^{a\tau^{*}} = 1/e$, and then decreases for further increase in $\tau$.}
\label{fig:roc}
\end{figure}

% {\color{red}Figure 1 goes here...} \\ \\ \\ \\ \\ \vspace{30mm}\\

% We now characterise the dependence of the convergence characteristics on coefficient $b$, i.e. $k$. For this analysis, we consider the coefficient a (i.e., ) and the time delay τ ≥ 0
% to be constant. The analysis is similar to the analysis
% described above. We state the following results for the
% dependence of the rate of convergence on coefficient b.

% \begin{equation}
%  u = (−a + \sigma)\tau \tan(u),
%  \end{equation}
\subsection{Non-oscillatory convergence}
In addition to a faster rate of convergence, it is also required to have the system equilibrates without oscillations. Non-oscillatory convergence is a desirable characteristic in the design of dynamical systems. In this sub-section, we derive a necessary and sufficient condition for non-oscillatory convergence of \eqref{eq:bd_delayedmodel}. For the system to be non-oscillatory, the eigenvalues should be negative real numbers. Therefore, we seek
conditions on system parameters for which the characteristic equation \eqref{eq:bd_chareq_roc} has negative real solution.

% Recall that the characteristic equation of the RCP without queue size feedback is
% \begin{equation}
% \label{eq:ce_withoutq1}
% \lambda + \left(\frac{a}{\tau}\right)e^{-\lambda \tau} = 0.
% \end{equation}
Substituting $\lambda=-\sigma+j\omega$ in \eqref{eq:bd_chareq_roc} gives 
\begin{align}
 \sigma&=a + b e^{\sigma\tau}\cos(\omega\tau), \label{eq:nocsigma}\\                         
 \omega&= b e^{\sigma\tau}\sin(\omega\tau). \label{eq:nocomega}
\end{align}
Solving the equations \eqref{eq:nocsigma} and \eqref{eq:nocomega} yields
\begin{equation}
 (\sigma-a) \tau=\frac{\omega\tau}{\tan(\omega\tau)}.\label{eq:nocnoq1}
\end{equation}
%begin
For the non-oscillatory convergence to the equilibrium, we require the eigenvalues to be negative real numbers, i.e., the real parts of all the roots of the characteristic equation should be negative, and the imaginary parts of all the roots of the characteristic equation should be zero ($\omega=0$). The right-hand side of \eqref{eq:nocnoq1} is a decreasing function of $\omega$, and has a maximum value of $1$ at $\omega=0$. Therefore, if the left-hand side of \eqref{eq:nocnoq1} is greater than or equal to 1, then $\omega = 0$ is the only solution to this equation. Hence, the necessary and sufficient condition for the imaginary parts of all the characteristic roots to be zero is $(\sigma - a)\tau \geq 1$.
Now, \eqref{eq:nocomega} can be rewritten as 
\begin{equation}
 b\tau e^{\sigma\tau} \frac{\sin(\omega\tau)}{\omega\tau}=1.
\end{equation}
Taking the limit $\omega\rightarrow 0$ gives 
\begin{equation}
 \lim_{\omega\rightarrow 0}\ b\tau e^{\sigma\tau} \frac{\sin(\omega\tau)}{\omega\tau}=b \tau e^{\sigma\tau},
\end{equation}
and hence $b \tau e^{\sigma\tau}=1$. Since $(\sigma - a) \tau\geq 1$, then $b\tau e ^{a \tau}\leq(1/e)$.
Thus, the necessary and sufficient condition for non-oscillatory convergence of \eqref{eq:bd_delayedmodel} is
\begin{equation}
b \tau e^{a \tau} \leq \frac{1}{e}.
\end{equation}
%begin
We can verify this condition using  \eqref{eq:nocsigma} as follows. 
Taking the limit $\omega\rightarrow 0$ in  \eqref{eq:nocsigma}  gives 
\begin{equation*}
\lim_{\omega\rightarrow 0}\ \cos(\omega\tau) \frac{b \tau e^{\sigma\tau}}{(\sigma - a)\tau}=\frac{b \tau e^{\sigma\tau}}{(\sigma - a)\tau} = 1.
\end{equation*}
We can rewrite the above equation as 
\begin{equation}
 (\sigma - a )\tau e^{-(\sigma- a )\tau}= b \tau e^{a \tau}. \label{eq:nonosc_final}
\end{equation}
The maximum value of the function $ f(x) = x e^{-x}$ is $1/e$. Therefore,  \eqref{eq:nonosc_final} has solutions if and only if $b \tau e^{a \tau} \leq 1/e$.

\subsection{Discussion}
After analyzing the condition for non-oscillatory convergence, and the rate of convergence, it can be deduced that the rate of convergence is maximum at the boundary of the non-oscillatory regime, i.e., $b\tau e^{a \tau}=1/e$. Thus, if we choose parameters such that  $b\tau e^{a \tau}=1/e$, then the system converges to the equilibrium quickly. The system becomes oscillatory when $b\tau e^{a\tau} > 1/e$, and hence the rate of convergence decreases. At $\tau = {\cos^{-1}(-a/b)}/{\sqrt{b^2 - a^2}}$, the system transits into an unstable regime, and hence the convergence rate would be zero. From the stability and convergence analyses, we can understand that there exist trade-offs between various system parameters. The stability study enables us to predict if the trajectories of the dynamical system would converge to the equilibrium, when subjected to a small perturbation. The results of convergence analysis enables us to tune the system parameters to make sure that the system converges quickly to a stable equilibrium.

% As can be observed in Figure \ref{fig:rocVsa} (a), for $a=0.1$ $(<1/e)$, the system shows over-damped behavior, i.e., reaches equilibrium without oscillating. At $a=1/e$, the system reaches equilibrium as quickly as possible without any oscillations. For $a \in(1/e, \pi/2)$, the system behaves in an under-damped manner, i.e., existence of convergent oscillations (see Figure \ref{fig:rocVsa} (b)). As shown in Figure \ref{fig:rocVsa} (c) and Figure \ref{fig:rocVsa} (d), the system loses stability at $a=\pi/2$, and the amplitude of the undamped oscillations increases as $a$ increases beyond $\pi/2$. 
The dependence of system behavior on various parameters is summarized in Table \ref{tab:roc}.
\begin{table}[hbtp]
\caption{Effect of the value of time delay ($\tau$) on the system behavior. The system reaches equilibrium quickly without any oscillations at $\tau = \tau^{*}$,  where $\tau^{*}$ is given by $b\tau^{*} e^{a \tau^{*}}=1/e$.}\label{tab:roc}
 % \begin{minipage}{\columnwidth}
\centering
\resizebox{0.95\columnwidth}{!}{%
\begin{tabular}{ll}
%  \toprule
\hline
    Parameter range&System behavior\vspace{1mm}\\
%    \midrule
\hline
    $\tau \in (0, \tau^{*}]$  & stable and non-oscillatory \vspace{1mm} \\ 
$\tau \in (\tau^{*}, \frac{\cos^{-1}(-a/b)}{\sqrt{b^2 - a^2}})$  & stable and oscillatory\vspace{1mm} \\
$ \tau \geq \frac{\cos^{-1}(-a/b)}{\sqrt{b^2 - a^2}}$  &  unstable\\
\hline
 % \bottomrule
  \end{tabular}
  }
%\end{center}
%\end{minipage}
\end{table}%

% {\color{red}To validate the theoretical results of non-oscillatory convergence, numerical simulations obtained using XPPAUT are shown in Figure 2. }
% 
% {\color{red}Figure 2 goes here..} \\ \vspace{50mm}\\
%    \begin{figure}[thpb]
%       \centering
%       %\includegraphics[scale=1.0]{figurefile}
%       \caption{Inductance of oscillation winding on amorphous
%        magnetic core versus DC bias magnetic field}
%       \label{figurelabel}
%    \end{figure}
So far, we have analyzed the some of the stability and convergence properties of \eqref{eq:bd_delayedmodel}. 

The linearized version of \eqref{eq:bd_delayedmodel_quadratic} is given by 
\begin{equation}
  \dot{u}(t) = \eta \big((-2 x_e + \mu) u(t)- k u(t-\tau)\big). \label{eq:bd_linearmodel_quadratic}
\end{equation}
From \eqref{eq:bd_linearmodel} and \eqref{eq:bd_linearmodel_quadratic}, we can note that the linearized model of \eqref{eq:bd_delayedmodel_quadratic} is quite similar to that of \eqref{eq:bd_delayedmodel}. Therefore, the results of linear analyses like stability, rate of convergence, and non-oscillatory convergence of \eqref{eq:bd_delayedmodel} can be extended to \eqref{eq:bd_delayedmodel_quadratic}. 

The next natural step is to investigate the dynamical behavior of the system as it transits from a stable to an unstable regime.  In the local stability analysis, we have shown that the system undergoes a Hopf bifurcation, as the bifurcation parameter crosses a critical value. In the next section, we study the characteristics of the bifurcating periodic solutions.
% So it is natural to study the characteristics of the bifurcating periodic solutions. To that end, we conduct a local Hopf bifurcation analysis in the next section.
\section{HOPF BIFURCATION}
% In the previous sections, we analyzed local stability and convergence properties of \eqref{eq:bd_delayedmodel}. It is also important to make sure that any loss of stability, that may happen, results in {stable} and small amplitude oscillations. 
% We hasten to add that we are not interested in destabilizing the system, but wish to employ the tools offered by local instability analysis to gain some insights into the non-linear properties of the system under consideration. 
In this section, we conduct a detailed Hopf bifurcation analysis for both the cubic and quadratic models.
% explore the impact of loss of local stability. In particular, we are concerned with the loss of local stability occurring via a Hopf bifurcation leading to the onset of limit cycles, as the bifurcation parameter crosses a critical value.
% In the bifurcation-theoretic analysis, we have to take non-linear terms into consideration, which helps to learn additional non-linear dynamical properties of the system.

Using the theoretical frameworks of Poincar\'{e} normal form and the center manifold theorem \cite{hassard1981} (which are outlined in the Appendix), we analytically characterize the type of the Hopf bifurcation and the stability of the bifurcating limit cycles.
% {\color{black}In addition, we also analyze the bifurcation properties of the quadratic model \eqref{eq:bd_delayedmodel_quadratic}, where the cubic term in \eqref{eq:bd_delayedmodel} is replaced by a quadratic term.}

The Hopf bifurcation analysis relies on both linear and non-linear terms of the Taylor series expansion of the non-linear model. As outlined in the Appendix, the stability and direction of the bifurcating limit cycles can be determined from the sign of first Lyapunov coefficient ($\mu_2$) and Floquet exponent ($\beta_2$), where
\begin{eqnarray}
\mu_2 &&\hspace{-6mm}= \dfrac{-\operatorname{Re}[c_1(0)]}{\alpha'(0)},\label{eq:alpha_2}\label{eq:alpa_2}\\
\beta_2 &&\hspace{-6mm}= 2\operatorname{Re}[c_1(0)].\label{eq:bta_2}
%\beta &&\hspace{-6mm}= \epsilon^2\beta_2+\mathcal{O}(\epsilon^4)\quad\hspace{-6mm} \beta_2 = 2\operatorname{Re}[c_1(0)]\quad \epsilon = \sqrt{\frac{\mu}{\mu_2}}\nonumber\\\label{49}
\end{eqnarray}
The sign of $\mu_2$ determines the direction of the Hopf bifurcation. If $\mu_2>0$ then the Hopf bifurcation is supercritical if $\mu_2<0$ it is subcritical. The sign of $\beta_2$ determines the stability of the bifurcating periodic solutions. The periodic solutions are asymptotically orbitally stable if $\beta_2 < 0$ and unstable if $\beta_2 >0$. Here, $\alpha'(0)$ is the real part of $d\lambda/d\eta$ evaluated at $\eta =  \eta_c$. We have already shown that $\alpha'(0) > 0$. Therefore, from \eqref{eq:alpa_2} and \eqref{eq:bta_2}, super-criticality of the bifurcating solution also establishes asymptotic orbital stability.

\subsection{Cubic model}
To recapitulate, the Taylor series expansion of \eqref{eq:bd_delayedmodel} is given by 
\begin{equation}
 \frac{d}{dt}u(t) = \eta \big( -a u(t) - b u(t - \tau) -3x_e u^2(t) - u^3(t)\big), \label{eq:bd_taylor_exp_recap}
\end{equation}
where $b = k $, $a = (3 x^2_e - \mu) $, and $x_e$ is given by $x_e^3 + (k- \mu) x_e + \Lambda= 0.$
Using the definitions outlined in the Appendix, the expression for $\mu_2$ of \eqref{eq:bd_delayedmodel} has been calculated as
 \begin{align}
\mu_2 = \frac{\xi^2_{xx}}{b^2}\ \tilde{g}(\epsilon)  + \frac{\xi_{xxx}}{b}\ \tilde{h}(\epsilon), \label{eq:bd_myu2_cubic_intermsofaandb}
\end{align}
where $\epsilon = \frac{a}{b}$ and 
 \begin{align}
\tilde{g}(\epsilon) = & \hspace{2mm} \dfrac{\sqrt{1 - \epsilon^2}(12\epsilon -18) + \cos^{-1}(-\epsilon) (8\epsilon^2 -18\epsilon +4)}{b^2 (1+\epsilon) (1-\epsilon^2) \cos^{-1}(-\epsilon) (5-4\epsilon)}
% &\Bigg( \dfrac{2\Big(1 + \frac{\epsilon \cos^{-1}(-\epsilon)}{\sqrt{1-\epsilon^2}} \Big) \big(1 - \epsilon - 2 \epsilon^2 \big) + 4(1 + \epsilon)\sqrt{1 - \epsilon^2} \cos^{-1}(-\epsilon)}{ {\big( 5 -3\epsilon^2 -4\epsilon^3 + 6\epsilon \big)}  ( \cos^{-1}(-\epsilon) \sqrt{1-\epsilon^2}) } \nonumber \\
% %  &+ \frac{2 \eta b \tau \sqrt{b^2 - a^2} \big(2 \sin(\theta) - \sin(2 \theta)\big)}{{\big(a + b ( 2 \cos(\theta) - \cos(2 \theta))\big)}^2 + \big(b (2 \sin(\theta) - \sin(2\theta)\big)^2} \nonumber \\
%  &- \frac{4 (1 + \frac{\epsilon \cos^{-1}(-\epsilon)}{\sqrt{1-\epsilon^2}})}{ (1 + \epsilon) \sqrt{1-\epsilon^2} \cos^{-1}(-\epsilon) }\Bigg) 
\label{eq:myu2_g}\\
 \tilde{h}(\epsilon) = & \hspace{2mm}  \frac{-3\sqrt{1 - \epsilon^2} -3 {\epsilon \cos^{-1}(-\epsilon)} }{ (1-\epsilon^2) \cos^{-1}(-\epsilon)} .\label{eq:myu2_h}
\end{align}

\begin{figure}[hbtp!]   
\vspace{-10mm}
\centering
%\mypsfrag{-0.6}{-0.2}\mypsfrag{0.0}{0.2}
\psfrag{0}{  \hspace{-2mm} $0$}
\psfrag{-50}{ \small \hspace{-2mm} $-50$}
\psfrag{-100}{ \small  \hspace{-3mm}  $-100$}
\psfrag{-150}{ \small \hspace{-2mm}  $-150$}
\psfrag{0.0}{ \small \hspace{-0mm} $0$}
\psfrag{0.2}{ \small \hspace{-1mm} $0.2$}
\psfrag{1.0}{ \small \hspace{-1mm} $1.0$}
\psfrag{0.4}{\small \hspace{-1mm} $0.4$}
\psfrag{0.6}{ \small \hspace{-1mm} $0.6$}
\psfrag{0.8}{ \small \hspace{-1mm} $0.8$}
\psfrag{pi/4}{ \small $\pi/4$}
\psfrag{pi/2}{ \small $\pi/2$}
\psfrag{x}{  $\epsilon$}
\psfrag{myu2}{  $\tilde{g}(\epsilon)$}
\psfrag{0.8}{ \small $0.8$}
%\psfrag{M}[b][t]{\hspace{0mm}\footnotesize $\mu_2$}
\psfrag{M}{\hspace{0mm}\footnotesize $\mu_{2}$}
\psfrag{B}{\hspace{0mm}\footnotesize $\beta_{2}$}
\includegraphics[height=2.5in,width=2.7in]{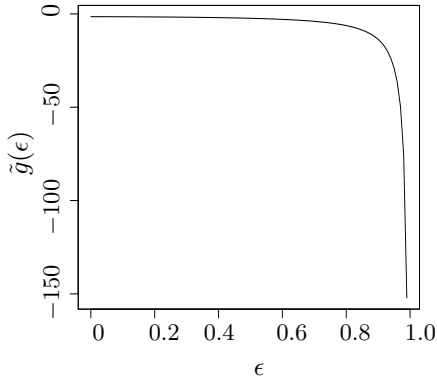}
\vspace{-3mm}
\caption{ Variation in $\tilde{g}(\epsilon)$ as $\epsilon$ is varied in the interval $[0,1)$. Observe that the sign of $\tilde{g}(\epsilon)$ is negative, and hence $\mu_2 < 0$, which implies that the Hopf bifurcation is sub-critical.}\label{fig:myu2plot_q}
\end{figure}

We now analyze the impact of quadratic and cubic terms on the type of the Hopf bifurcation.

As $b>a$, $a \geq 0$, and $b > 0$, the value of $\epsilon = \frac{a}{b}$ would lie in the interval $[0,1)$. From \eqref{eq:myu2_h}, we can readily show that $\tilde{h}(\epsilon) < 0$ and, hence $\xi_{xxx} \tilde{h}(\epsilon) > 0$, as $\xi_{xxx} = -1$. Therefore, the presence of cubic term induces super-criticality. To analyze the effect of quadratic term ($\xi_{xx}$) on the nature of the bifurcation, we plot the value of $\tilde{g}(\epsilon)$, as $\epsilon$ varies in the interval $[0,1)$. From Fig. \ref{fig:myu2plot_q}, we can observe that the value of $\mu_2 < 0$ for $\epsilon \in [0,1)$. Hence, the type of Hopf bifurcation is sub-critical.
% Substituting the values of the Taylor series coefficients $\xi_{xx}$ and $\xi_{xxx}$ in \eqref{eq:bd_myu2_cubic}, we get
% \begin{footnotesize}
% \begin{align}
%   \mu_2 = 9 x^{2}_e &\Bigg( \dfrac{2\big(1 + \eta a \tau\big) \big(a + b(2\cos(\theta) - \cos({2\theta})) \big)}{ {\big(a + b ( 2 \cos(\theta) - \cos(2 \theta))\big)}^2 + (b (2 \sin(\theta) - \sin(2\theta))^2 } \nonumber \\
%  &+ \frac{2 \eta b \tau \sqrt{b^2 - a^2} \big(2 \sin(\theta) - \sin(2 \theta)\big)}{{\big(a + b ( 2 \cos(\theta) - \cos(2 \theta))\big)}^2 + \big(b (2 \sin(\theta) - \sin(2\theta)\big)^2} \nonumber \\
%  &- \frac{4 (1 + \eta a \tau)}{ (a + b) \tau ( b^2 - a^2) }\Bigg)  + \Bigg( \frac{3(1 + \eta a \tau)}{\tau (b^2 - a^2)} \Bigg), \label{eq:bd_myu2_cubic_1}
% \end{align}
% \end{footnotesize}
% where  $a = (3 x^2_e - \mu)$,  $b = k$, and $x_e$ is given by $x_e^3 + (k- \mu) x_e + \Lambda= 0.$
%  \frac{d}{dt}u(t) = \eta \big((-3 x^2_e + \mu) u(t) - k u(t - \tau) -3x_e u^2(t) - u^3(t)\big). \label{eq:bd_taylor_exp}
% \end{equation}
%where $x_e$ is given by $x_e^3 + (k- \mu) x_e + \Lambda= 0.$
Therefore, when we have both $\xi_{xx}$ and $\xi_{xxx}$, then the type of the Hopf bifurcation depends on the values of the model parameters. We now validate the analytical results using some numerical examples where we choose specific values for the system parameters. \\

\begin{figure}[hbtp!]
\vspace{-8mm}
\centering
\psfrag{R}[][][2]{\small Amplitude of oscillation}
\psfrag{kappa}[][][2]{\small Bifurcation parameter, $\eta$}
\psfrag{0.95}[][][2]{ \footnotesize  $0.95$}
\psfrag{1.00}[][][2]{  \footnotesize  $1.00$}
\psfrag{1.01}[][][2]{  \footnotesize  $1.01$}
\psfrag{1.02}[][][2]{  \footnotesize  $1.02$}
\psfrag{1.03}[][][2]{  \footnotesize  $1.03$}
\psfrag{1.04}[][][2]{  \footnotesize  $1.04$}
\psfrag{1.05}[][][2]{  \footnotesize  $1.05$}
\psfrag{0}[][][2]{  \footnotesize  $0$}
\psfrag{3}[][][2]{  \footnotesize  $3$}
\psfrag{6}[][][2]{  \footnotesize  $6$}
\includegraphics[scale = 0.9,trim=0cm 0cm 0cm 1.7cm, clip=true,width=3.25in,height=1.75in]{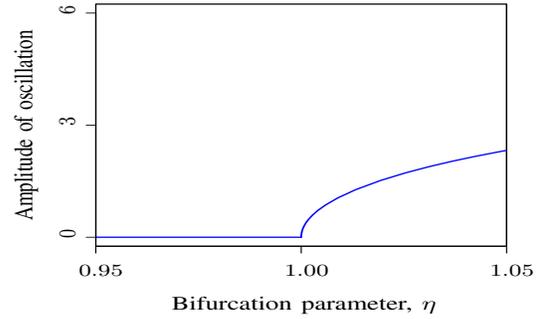}
\vspace{-1mm}
\caption{Bifurcation diagram highlighting that the system undergoes a super-critical Hopf bifurcation at $\eta=1$. The parameter values used are $k = 9$, $\tau = 0.187$, $\Lambda = -7$, and $\mu = 1$.}
%\vspace{-5mm}
\label{fig:chemkinetics_bfd_super}
\vspace{-2mm}
\end{figure}

$\textit{Numerical Example 1 (Super-critical)}$: Let us consider the system with $k = 9$, $\tau = 0.187$, $\Lambda = -7$, and $\mu = 1$. For these values, using the Hopf condition, we get $\eta_c =1$. From \eqref{eq:bd_myu2_cubic_intermsofaandb}, we obtain the value of $\mu$ as $9.3$, which implies that the system undergoes a super-critical Hopf bifurcation. The bifurcation diagram drawn using the Matlab package DDE-Biftool \cite{ddetool1,ddetool2} is shown in Figure \ref{fig:chemkinetics_bfd_super}. As expected, it shows that the system loses local stability via a super-critical Hopf bifurcation, as the bifurcation parameter crosses the critical threshold ($\eta_c =1$). To validate this, numerical simulations obtained using XPPAUT \cite{xppaut2002} are shown in Figure \ref{fig:chemkineticstsplots_super12}. For $\eta=0.95$, the system converges to the equilibrium (see Figure \ref{fig:chemkineticstsplots_super12}(a)). Whereas, for $\eta = 1.05 > \eta_c$ i.e. after the bifurcation, the system leads to the emergence of stable and small-amplitude limit cycles.\\

\newcommand{\supwdth}{0.25\textwidth}
\begin{figure}[hbtp!]
\vspace{-8mm}
\psfrag{time}{\hspace{-0.1cm}  \footnotesize Time}
\psfrag{Rate}{\hspace{-0.53017cm} $x(t)$}
\psfrag{0}{\footnotesize{$0$}}
\psfrag{5}{\hspace{-0.05cm}\scriptsize{$5$}}
\psfrag{10}{\hspace{-0.1cm}\scriptsize{$10$}}
\psfrag{15}{\hspace{-0.1cm}\scriptsize{$15$}}
\psfrag{20}{\hspace{-0.1cm}\scriptsize{$20$}}
\psfrag{25}{\hspace{-0.1cm}\scriptsize{$25$}}
\psfrag{50}{\hspace{-0.1cm}\scriptsize{$50$}}
\psfrag{10000}{\hspace{-0.1cm}\scriptsize{$10000$}}
\psfrag{20000}{\hspace{-0.1cm}\scriptsize{$20000$}}
\psfrag{0.5}{\hspace{-0.1cm}\scriptsize{$0.95$}}
\psfrag{1.0}{\hspace{-0.1cm}\scriptsize{$1.0$}}
\psfrag{1.5}{\hspace{-0.1cm}\scriptsize{$1.05$}}
\psfrag{3}{\hspace{-0.1cm}\scriptsize{$3$}}
\psfrag{-3}{\hspace{-0.1cm}\scriptsize{$-3$}}
\psfrag{1}{\hspace{-0.1cm}\scriptsize{$1$}}
\psfrag{2}{\hspace{-0.1cm}\scriptsize{$2$}}
\begin{tabular}{c}
\subfloat[  $\eta = 0.95$]{\includegraphics[width=\supwdth]{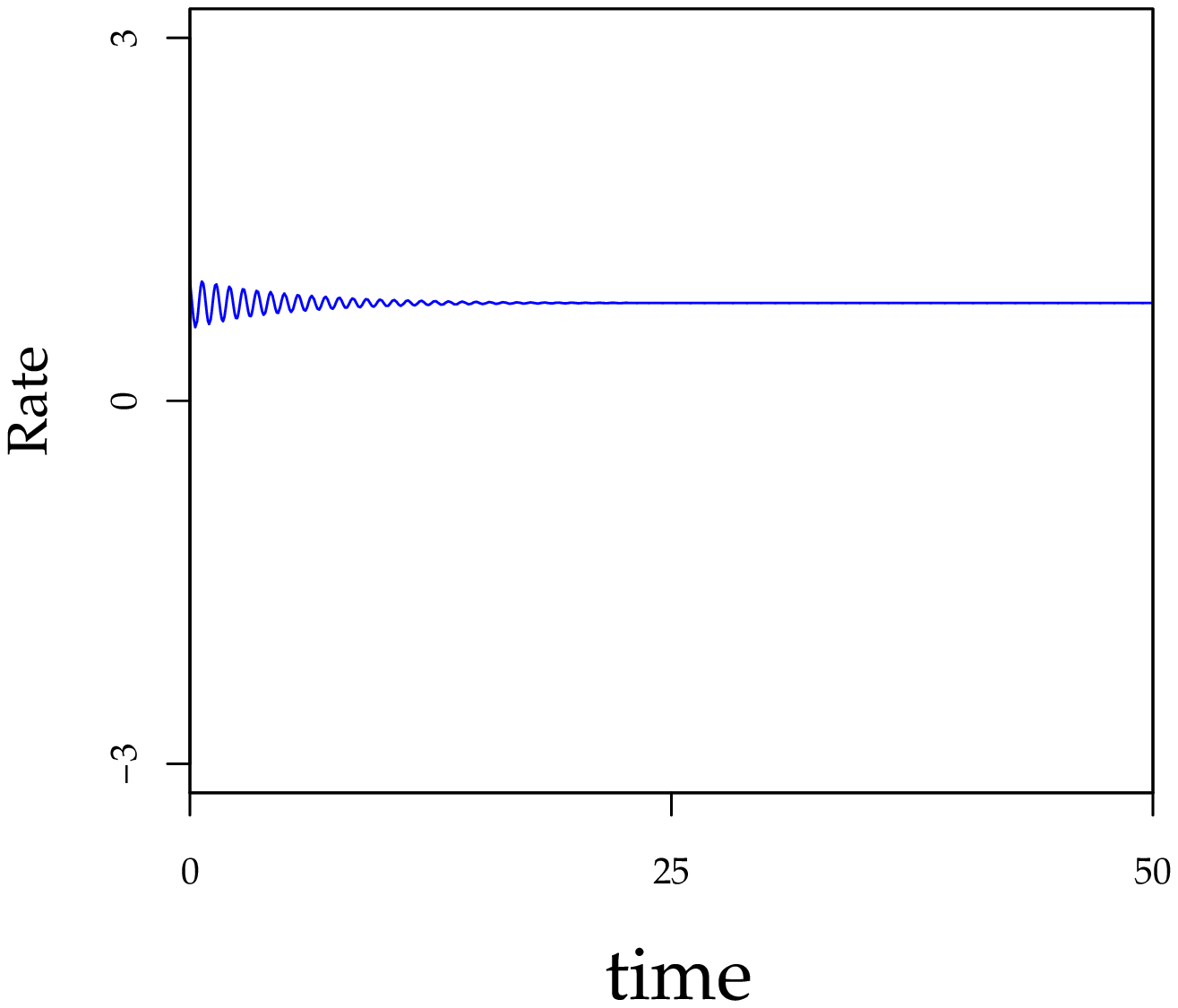}}\hspace{-5mm}
%\subfloat[ $\eta = 0.97, R_0=2$]{\includegraphics[width=\supwdth]{tompecssuper2.eps}}\hspace{-5mm}
\subfloat[ $\eta = 1.05$]{\includegraphics[width=\supwdth]{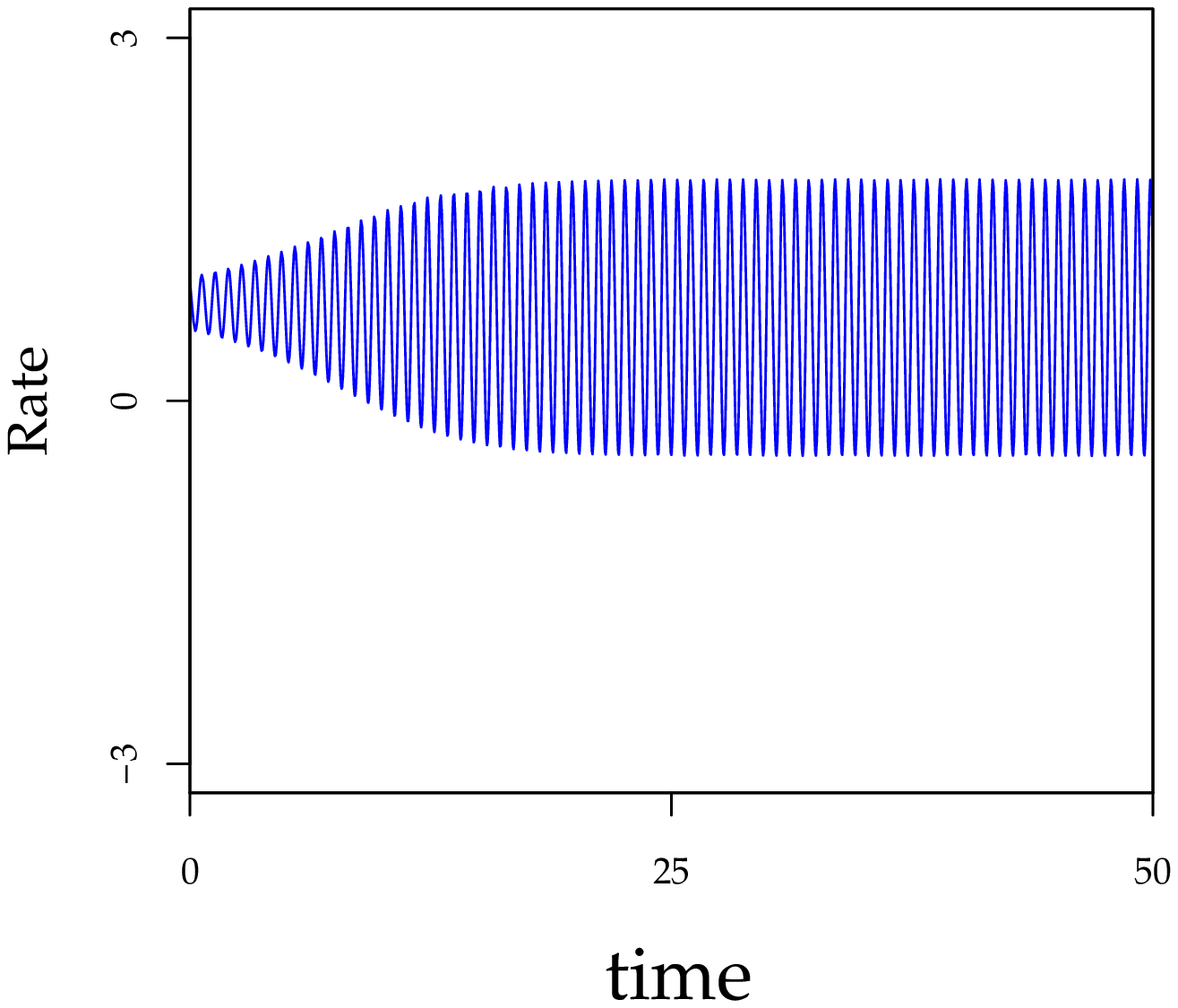}}
\end{tabular}
\caption{Numerical computations illustrating that the system exhibits a super-critical Hopf bifurcation, as $\eta$ increases  beyond the critical value. Time series plots are shown for the cases $\eta < 1$ and $\eta > 1$. The parameter values chosen are $k = 9$, $\tau = 0.187$, $\Lambda = -7$, and $\mu = 1$. Here, we set initial condition as $x_0 = 0.9.$}
\label{fig:chemkineticstsplots_super12}
\end{figure}

% We now resort to numerical bifurcation analysis tools. We analyze the nature of the Hopf bifurcation by plotting the bifurcation diagrams numerically (using DDE-Biftool), and then validate it through numerical simulations (using XPPAUT \cite{xppaut2002}). We consider some special cases where we choose specific values for the system parameters.  
% We now analyze what happens when we have both $\xi_{xx}$ and $\xi_{xxx}$. Substituting the values of $\xi_{xx}$ and $\xi_{xxx}$ in \eqref{eq:bd_myu2_cubic_intermsofaandb}, we get
% $\textit{Numerical Example 2 (Sub-critical)}$: Let $a=0.924$, $C=10$, $\tau=100$ and $b=0.257$ which corresponds to equilibrium utilization of 70\% ($\rho^*=0.70$), the system undergoes a Hopf bifurcation at $\eta$ = 1. Using \eqref{eq:myu2rcp2}, we calculate $\mu_2  = -3.547 \times 10^{-2}  <0$, implying that the Hopf bifurcation is sub-critical. The bifurcation diagram shown in Fig. \ref{fig:subcriticalplot} confirms that the system exhibits a sub-critical Hopf bifurcation, as the bifurcation parameter is varied beyond the critical threshold.
\begin{figure}[hbtp!]
\vspace{-6mm}
\centering
\psfrag{R}[][][2]{\small Amplitude of oscillation}
\psfrag{kappa}[][][2]{\small Bifurcation parameter, $\eta$}
\psfrag{0.95}[][][2]{ \footnotesize $0.95$}
\psfrag{1.00}[][][2]{ \footnotesize $1.00$}
\psfrag{1.01}[][][2]{ \footnotesize  $1.01$}
\psfrag{1.02}[][][2]{ \footnotesize  $1.02$}
\psfrag{1.03}[][][2]{  \footnotesize $1.03$}
\psfrag{1.04}[][][2]{  \footnotesize $1.04$}
\psfrag{1.05}[][][2]{  \footnotesize $1.05$}
\psfrag{0.70}[][][2]{ \footnotesize $0.70$}
\psfrag{0.80}[][][2]{  \footnotesize $0.80$}
\psfrag{0.90}[][][2]{  \footnotesize $0.90$}
\psfrag{1.0}[][][2]{  \footnotesize $1.0$}
\psfrag{1.1}[][][2]{  \footnotesize $1.1$}
\psfrag{0}[][][2]{  \footnotesize $0$}
\psfrag{3}[][][2]{  \footnotesize $3$}
\psfrag{6}[][][2]{  \footnotesize $6$}
\includegraphics[scale = 0.9,trim=0cm 0cm 0cm 1.7cm, clip=true,width=3.25in,height=1.75in]{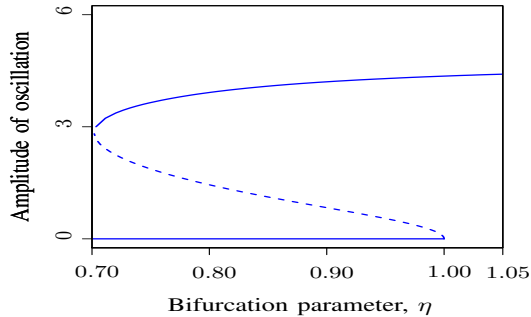}
\vspace{-1mm}
\caption{Bifurcation diagram showing the existence of a sub-critical Hopf for the parameter values $k = 4.75$, $\tau = 1$, $\mu = 1$, and $\Lambda= -7$. The solid and dashed lines denote the amplitude of stable and unstable limit cycles, respectively. }
\label{fig:chemkinetics_bfd_sub}
\end{figure}
\begin{figure}[hbtp!]
\newcommand{\subwdth}{0.25\textwidth}
\psfrag{time}{\hspace{-0.1cm}  \footnotesize Time}
\psfrag{Rate}{\hspace{-0.57cm} $x(t)$}
\psfrag{0}{\hspace{-0.04cm}\scriptsize{$0$}}
\psfrag{5}{\hspace{-0.05cm}\scriptsize{$5$}}
\psfrag{10}{\hspace{-0.1cm}\scriptsize{$10$}}
\psfrag{15}{\hspace{-0.1cm}\scriptsize{$15$}}
\psfrag{20}{\hspace{-0.1cm}\scriptsize{$20$}}
\psfrag{25}{\hspace{0cm}\scriptsize{$25$}}
\psfrag{50}{\hspace{-0.1cm}\scriptsize{$50$}}
\psfrag{100}{\hspace{-0.1cm}\scriptsize{$100$}}
\psfrag{200}{\hspace{-0.1cm}\scriptsize{$200$}}
\psfrag{-3}{\hspace{-0.1cm}\scriptsize{$-3$}}
\psfrag{3}{\hspace{-0.1cm}\scriptsize{$3$}}
\psfrag{1.5}{\hspace{-0.1cm}\scriptsize{$1.05$}}
\psfrag{150}{\hspace{-0.1cm}\scriptsize{$150$}}
\psfrag{300}{\hspace{-0.1cm}\scriptsize{$300$}}
\psfrag{400}{\hspace{-0.1cm}\scriptsize{$400$}}
\psfrag{125}{\hspace{-0.1cm}\scriptsize{$125$}}
\psfrag{250}{\hspace{-0.1cm}\scriptsize{$250$}}
\psfrag{60}{\hspace{-0.1cm}\scriptsize{$60$}}
\psfrag{120}{\hspace{-0.1cm}\scriptsize{$120$}}
\begin{tabular}{c}
\subfloat[ $\eta = 0.95, x_0 = 3$]{\includegraphics[width=\subwdth]{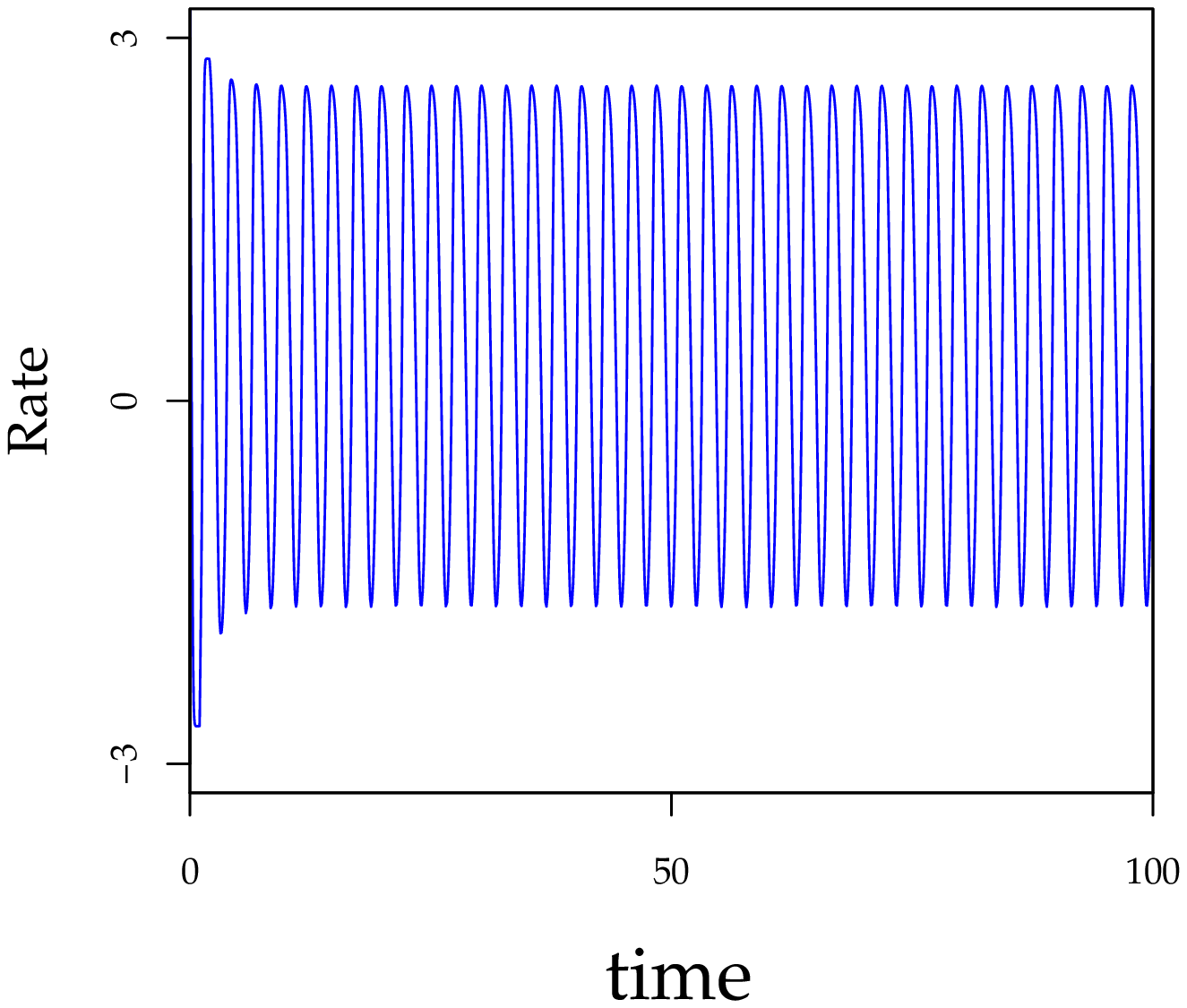}}\hspace{-5mm}
% \subfloat[ $\eta = 0.97, R_0=2$]{\includegraphics[width=\subwdth]{tompecssub2.eps}}\hspace{-5mm}
\subfloat[ $\eta = 1.05, x_0 = 1.35$]{\includegraphics[width=\subwdth]{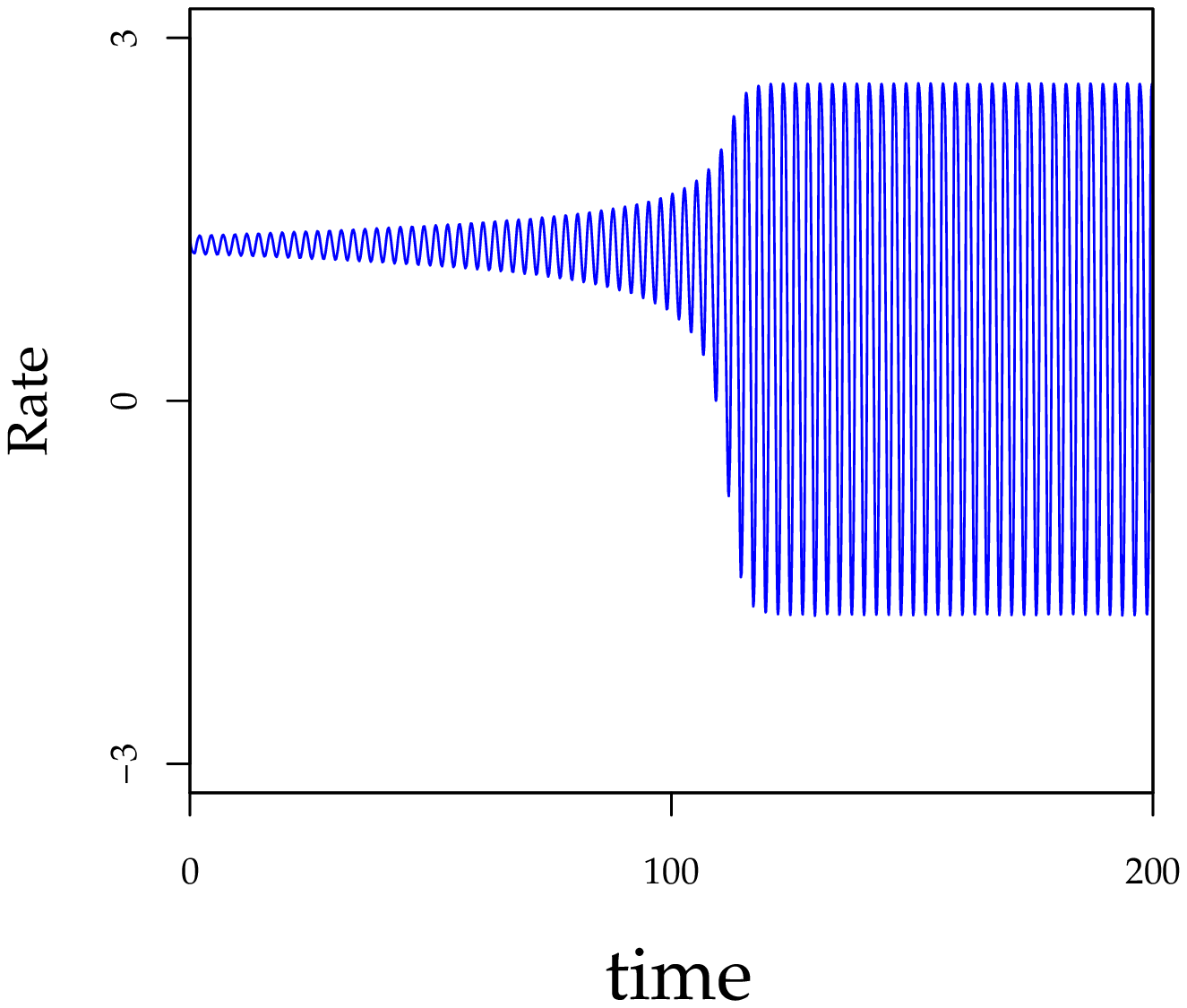}}
\end{tabular}
\caption{Numerical computations illustrating that the system undergoes a sub-critical Hopf as the bifurcation parameter $\eta$ increases beyond the critical threshold ($\eta_c =1$). The values of the parameters used are $k = 4.75$, $\tau = 1$, $\mu = 1$, and $\Lambda= -7$. For these values, we obtain equilibrium point as $x_e = 1.3$.} \label{fig:chemkineticstsplots_sub}
\end{figure}

% The numerical computation shown in Fig. \ref{fig:sub_rcp1d_withq}(a) illustrates that the system is stable for $\eta<1$. For $\eta = 1.05$, as shown in Fig. \ref{fig:sub_rcp1d_withq}(b), the system exhibits a limit cycle but now with amplitude larger than that of previous example.
$\textit{Numerical Example 2 (Sub-critical)}$: Consider $k = 4.75$, $\tau = 1$, $\mu = 1$, and $\Lambda= -7$. The system undergoes a Hopf bifurcation at $\eta$ = 1. For these values, we get $\mu_2 = -6.11 < 0$, implying that the Hopf bifurcation is sub-critical. The bifurcation diagram shown in \ref{fig:chemkinetics_bfd_sub} confirms that the system exhibits a sub-critical Hopf bifurcation, as the bifurcation parameter is varied beyond the critical threshold. 
The numerical computation shown in Figure \ref{fig:chemkineticstsplots_sub} illustrates that the system is stable for $\eta < 1$. For $\eta = 1.05$, as shown in Figure \ref{fig:chemkineticstsplots_sub}, the system exhibits a limit cycle but now with amplitude larger than that of previous example.\\
%%%%%%%%%%%%%%%%%%%%%%%%%%%%%%%%%%%%%%%%%%%%%%%%%%%%%

\subsection{Quadratic model}
The Taylor series expansion of the quadratic model \eqref{eq:bd_delayedmodel_quadratic} about its equilibrium is given by 
\begin{equation}
  \frac{d}{dt}u(t) = \eta \big(-a u(t)- b u(t-\tau) - u^2(t) \big), \label{eq:bd_taylor_exp_quadratic}
\end{equation}
where $a = (2 x_e - \mu)$, $b = k$, and $x_e$ is given by $x_e^2 + (k- \mu) x_e + \Lambda= 0.$ For this quadratic model, we obtain $\mu_2$ as 
% 
% %{\color{black}
% In the case of quadratic model \eqref{eq:bd_delayedmodel_quadratic}, the value of $\mu_2$ is given by
% \begin{small}
%  \begin{equation}
%  \mu_2 = \xi^2_{xx} \Bigg( \dfrac{\sqrt{1 - \epsilon^2}(12\epsilon -18) + \cos^{-1}(-\epsilon) (8\epsilon^2 -18\epsilon +4)}{b^2 (1+\epsilon) (1-\epsilon^2) \cos^{-1}(-\epsilon) (5-4\epsilon)} \Bigg)\\ 
% %  \Bigg( \dfrac{2\big(1 + \eta a \tau\big) \big(a + b(2\cos(\theta) - \cos({2\theta})) \big)}{ {\big(a + b ( 2 \cos(\theta) - \cos(2 \theta))\big)}^2 + (b (2 \sin(\theta) - \sin(2\theta))^2 } \nonumber \\
% %  &+ \frac{2 \eta b \tau \sqrt{b^2 - a^2} \big(2 \sin(\theta) - \sin(2 \theta)\big)}{{\big(a + b ( 2 \cos(\theta) - \cos(2 \theta))\big)}^2 + \big(b (2 \sin(\theta) - \sin(2\theta)\big)^2} \nonumber \\
% %  &- \frac{4 (1 + \eta a \tau)}{ (a + b) \tau ( b^2 - a^2) }\Bigg),
%  %+ \xi_{xxx} \Bigg( \frac{-3(1 + \eta a \tau)}{\tau (b^2 - a^2)} \Bigg), 
%  \label{eq:bd_myu2_quadratic}
% \end{equation}
% \end{small}
% where $\theta = \cos^{-1}(-a/b)$. This can be further simplified as 
%\begin{footnotesize}
 \begin{align}
\mu_2 = \frac{\xi^2_{xx}}{b^2}\ \tilde{g}(\epsilon) = \frac{\tilde{g}(\epsilon)}{b^2}, \label{eq:bd_myu2_quad_intermsofaandb}
\end{align}
%\end{footnotesize}
where $\tilde{g}(\epsilon)$ is given by \eqref{eq:myu2_g}. We have already shown that $\tilde{g}(\epsilon) < 0$ for all values $\epsilon = a/b \in [0,1)$. Therefore, the quadratic model would undergo a sub-critical Hopf bifurcation, which would result in the emergence of either large amplitude limit cycles or unstable limit cycles. Both these outcomes are undesirable, and hence the sub-critical Hopf should be avoided.

In general, it would be preferable to have a stable equilibrium. However, if the system does lose stability due to variation in system parameters, it would be desirable to have an asymptotically orbitally stable limit cycle of small amplitude. To that end, a super-critical Hopf may be preferable over a sub-critical Hopf bifurcation. The Boissonade-De-Kepper model \eqref{eq:bd_delayedmodel}, which incorporates cubic control law, can undergo both super-critical and sub-critical Hopf bifurcation, depending on the parameter values. Whereas, in the case of quadratic model \eqref{eq:bd_delayedmodel_quadratic}, the Hopf bifurcation is sub-critical. Thus, our results tend to favor the cubic model. 

% Substituting the value of $\xi_{xx}$ in \eqref{eq:bd_myu2_cubic}, we get
% \begin{footnotesize}
% \begin{align}
%   \mu_2 = &\Bigg( \dfrac{2\big(1 + \eta a \tau\big) \big(a + b(2\cos(\theta) - \cos({2\theta})) \big)}{ {\big(a + b ( 2 \cos(\theta) - \cos(2 \theta))\big)}^2 + (b (2 \sin(\theta) - \sin(2\theta))^2 } \nonumber \\
%  &+ \frac{2 \eta b \tau \sqrt{b^2 - a^2} \big(2 \sin(\theta) - \sin(2 \theta)\big)}{{\big(a + b ( 2 \cos(\theta) - \cos(2 \theta))\big)}^2 + \big(b (2 \sin(\theta) - \sin(2\theta))\big)^2} \nonumber \\
%  &- \frac{4 (1 + \eta a \tau)}{ (a + b) \tau ( b^2 - a^2) }\Bigg). \label{eq:bd_myu2_quadratic_1}
% \end{align}
% \end{footnotesize}
% }

%%%%%%%%%%%%%%%%%%%%%%%%%%%%%%%%%%%%%%%%%%%%%%%%%%%%%%%%%%%%%%%%%%%%%%%%%%%%%%%%
\section{CONTRIBUTIONS}
The contribution of this paper is two fold. First, mathematicians concerned with the investigation of delay differential equations would be given a stimulating example of application in chemistry; second, chemists would further realize that on all levels of chemical processes feedback delays play an important role in the generation of instabilities.

We analyzed the stability, convergence, and bifurcation properties of the Boissonade-De Kepper (BD) model with feedback delay. 
% , which provides the theoretical underpinning for the systematic search procedure that has led to the discovery of many new oscillatory reactions. 
From the results of stability and convergence analyses, one can tune the system parameters to make sure that the system converges quickly to a stable equilibrium. We also showed that the system  undergoes a Hopf bifurcation, if the stability condition gets violated. Using Poincar\'{e} normal form and the center manifold theorem \cite{hassard1981}, we analyzed the direction and stability of the bifurcating limit cycles. We also investigated the bifurcation properties of a quadratic model, where the cubic term in the original BD model is replaced by a quadratic term. We established that the BD model, which incorporates cubic control law, can undergo both super-critical and sub-critical Hopf, depending on the parameter values. Whereas, in the case of quadratic model, the Hopf bifurcation is always sub-critical, which may give rise to either limit cycles with a large amplitude or unstable limit cycles. Therefore, our results tend to favor the cubic model. We validated some of our analytical insights using  bifurcation diagrams and numerical simulations.

We also derived general results useful in the study of the nature of the Hopf bifurcation of a general first-order non-linear delay differential equation. Therefore, the bifurcation results are not just confined to the BD model, but can also be extended to other non-linear delayed systems as well. To highlight the implications of our Hopf bifurcation results, we also applied our results to determine the nature of Hopf bifurcation of the Nicholson's Blowflies equation \cite{nicholsonbfref},  
% \subsection{Conclusions}
% 
% \subsection{Future Works}

%\section{APPENDIX}
\appendix
%\appendices

%\section{Hopf bifurcation analysis}
% % you can choose not to have a title for an appendix
% % if you want by leaving the argument blank
% \section{}

Here, we outline the necessary calculations to determine the type of Hopf bifurcation and the asymptotic form of the bifurcation solutions as local instability just sets in. For now, we will only be concerned with the first Hopf bifurcation.
The framework employed to address the stability of the limit cycles is the Poincar{\'e} normal form, and the center manifold theorem. 

Consider the following non-linear delay differential equation:
\begin{equation}
\label{eq:gen_non_eq}
\dfrac{d}{dt}x(t) = \eta f\big(x(t),x(t-\tau)\big),
\end{equation}
where $f$ has a unique equilibrium denoted by $x^*$ and $\tau > 0$. Define $u(t) = x(t)-x^*,$ and take a Taylor expansion for (\ref{eq:gen_non_eq}) including the linear, quadratic and cubic terms to obtain
\vspace{-1.5mm}
% \begin{equation}
% \begin{aligned}
% \dfrac{d}{dt}u(t)&= \eta \big(\xi_yu(t-\tau)+\xi_{xy}u(t)u(t-\tau)\\
% %&+\ \xi_{xy}u(t)u(t-\tau)\\
% &+\ \xi_{yy}u^2(t-\tau)\\
% &+\ \xi_{xyy}u(t)u^2(t-\tau)\\
% &+\ \xi_{yyy}u^3(t-\tau)\\
% &+\ \xi_{yyz}u^2(t-\tau)u(t-\tau_2)+\xi_{zzz}u^3(t-\tau_2)\\
% &+\ \xi_{yzz}u(t-\tau)u^2(t-\tau_2)+ \mathcal{O}(u^4)\big)\\
% \end{aligned}
% \label{eq:linear_noneq}
% \end{equation}

\begin{align}
 \frac{d}{dt}u(t) =&\,\, \eta(\xi_x u(t) + \xi_y u(t-\tau) + \xi_{xx}u^2(t)  \notag \\
 &+ \xi_{xy}u(t)u(t-\tau)+ \xi_{yy}u^2(t-\tau) + \xi_{xxx}u^3(t) \notag  \\
		  & + \xi_{xxy}u^2(t)u(t-\tau)  + \xi_{xyy}u(t)u^2(t-\tau) \notag  \\
		   &+ \xi_{yyy}u^3(t-\tau) ),
\label{eq:linear_noneq}
		  \end{align}

where, letting $f^*$ denote evaluation of $f$ at $(x^*,y^*)$
\begin{alignat*}
\xi\xi_i&=f^*_i,&\xi_{ii}&=\dfrac{1}{2}f^*_{ii},&\xi_{iii}&=\dfrac{1}{6}f^*_{iii} \quad \forall \  i \in \{x,y\}\\
\xi_{xy}&=f^*_{xy},\ & \xi_{xyy}&=\dfrac{1}{2}f^*_{xyy},\ & \xi_{xxy}&=\dfrac{1}{2}f^*_{xxy}.
\end{alignat*}
Considering the linearised form of (\ref{eq:linear_noneq}), we get
\begin{equation}
 \frac{d}{dt}u(t) = \eta\xi_{x}u(t) + \eta \xi_{y}u(t-\tau).
\label{eq:leqn}
  \end{equation}
Looking for exponential solutions, the characteristic equation of (\ref{eq:leqn}) is given by
\begin{equation}
\lambda + \eta a+\eta b e^{-\lambda\tau} = 0,
\end{equation}
where $a = -\xi_x$ and $b = -\xi_y$. We assume that $a \geq 0$, $ b > 0$, and $b > a$. Let $\eta=\eta_c+\mu$, and the system undergoes Hopf bifurcation at $\mu=0$.
The calculations that follow will enable us to address questions about the form of the bifurcating solutions, as the system transits
from stability to instability via a Hopf bifurcation. For this we have to take higher order terms, i.e., the quadratic and cubic of (\ref{eq:linear_noneq}) into consideration. Following the work of  \cite{hassard1981}, we now perform the requisite calculations.

Consider the following autonomous delay-differential system
\begin{equation}
\label{eq:auto_noneq}
\dfrac{d}{dt}u(t) = \mathcal{L}_\mu u_t + \mathcal{F}(u_t,\mu),
\end{equation}
where $t>0,\ \mu \in \mathbb{R},\ \tau >0,$
\begin{equation*}
  u_t(\theta) = u(t+\theta),\  u:[-\tau,0]\rightarrow\mathbb{R},\  \theta\in[-\tau,0].
\end{equation*}
$\mathcal{L}_\mu$ is a one-parameter family of continuous linear operators defined as $\mathcal{L}_\mu: C[-\tau,0]\rightarrow\mathbb{R}$. The operator $\mathcal{F}(u_t,\mu):C[-\tau,0]\rightarrow\mathbb{R}$  contains the non-linear terms. Further, assume that $\mathcal{F}(u_t,\mu)$ is analytic and that $\mathcal{F}$ and $\mathcal{L}_\mu$ depend analytically on the bifurcation parameter. Note that (\ref{eq:linear_noneq}) is a
type of the form (\ref{eq:auto_noneq}). The objective now is to rewrite (\ref{eq:auto_noneq}) as follow
\begin{equation}
\label{eq:auto_noneq_matrix}
 \dfrac{d}{dt}u_t  = \mathcal{A}(\mu)u_t+\mathcal{R}u_t
\end{equation}
which has $u_t$ rather than both $u$ and $u_t$. By the Riesz representation theorem, there exists a matrix-valued function $\zeta(.,\mu):[-\tau,0]\rightarrow\mathbb{R}^{n^2}$, with variation of each component of $\zeta$ is bounded and for all $\phi \in C[-\tau,0]$

\begin{equation*}
  \mathcal{L}_\mu\phi = \int_{-\tau}^0d\zeta(\theta,\mu)\phi(\theta),
\end{equation*}
% In particular
% \begin{equation}
% \label{eq:auto_noneq1}
%   \mathcal{L}_\mu u_t = \int_{-\tau}^0d\zeta(\theta,\mu)u(t+\theta).
% \end{equation}
where   $\mathrm{d}\zeta(\theta,\mu) =\eta\big(\xi_y\delta(\theta+\tau)+\xi_z\delta(\theta+\tau_2)\big)\mathrm{d}\theta$,
% \begin{equation*}
%   d\zeta(\theta,\mu) =\eta\big(\xi_y\delta(\theta+\tau)+\xi_z\delta(\theta+\tau_2)\big)d\theta,
% \end{equation*}
and $\delta(\theta)$ is the Dirac delta function. Now we define
\begin{equation}
\label{eq:Atheta}
  \mathcal{A}(\mu)\phi(\theta) =
  \begin{cases}
  \dfrac{d\phi(\theta)}{d\theta}, & \theta \in [-\tau,0),\\
  \int_{-\tau}^0d\zeta(s,\mu)\phi(s), & \theta=0,
  \end{cases}
\end{equation}
and
\begin{equation*}
  \mathcal{R}\phi(\theta) =
    \begin{cases}
  0, & \theta \in [-\tau,0),\\
  \mathcal{F}(\phi,\mu),& \theta=0.
  \end{cases}
\end{equation*}
Now the system (\ref{eq:auto_noneq}) becomes equivalent to (\ref{eq:auto_noneq_matrix}) as required.

% The bifurcating periodic solutions $u(t,\mu(\epsilon))$ of (\ref{eq:auto_noneq}) (where $\epsilon \geq 0$
% is a small parameter) have amplitude $\mathcal{O}(\epsilon)$, period $\mathcal{P}(\epsilon)$
% and nonzero Floquet exponent $\beta(\epsilon)$, where $\mu,\mathcal{P}$ and $\beta$ have the
% following (convergent) expansions:
% \begin{equation*}
% \begin{aligned}
%   &\mu = \mu_2\epsilon^2+ \mu_4\epsilon^4+\cdots\\
%   &\mathcal{P} = \dfrac{2\pi}{\omega_0}(1+\mathcal{T}_2\epsilon^2+\mathcal{T}_4\epsilon^4+\cdots) \\
%   &\beta = \beta_2\epsilon^2+ \beta_4\epsilon^4+\cdots.\\
% \end{aligned}
% \end{equation*}
% The sign of $\mu_2$ determines the \emph{direction of bifurcation}: If $\mu_2>0$, the bifurcation is \emph{supercritical} and $\mu_2<0$ implies a \emph{subcritical} bifurcation. The sign of $\beta_2$ determines the stability of $u(t,\mu(\epsilon))$: \emph{asymptotic orbital stability} if $\beta_2<0$ and instability if $\beta_2>0$. These coefficients will now be determined. 
%We only need to compute the expressions at $\mu=0$, hence, we set $\mu=0$ in the following. 
Let $q(\theta)$ be the eigenfunction for $\mathcal{A}(0)$ corresponding to $\lambda(0)$, namely
% \begin{equation*}
$\mathcal{A}(0)q(\theta) = i\omega_0q(\theta).$
% \end{equation*}
Now we define an adjoint operator $\mathcal{A}^*(0)$ as
\begin{equation*}
  \mathcal{A}^*(0)\alpha(s) =
    \begin{cases}
  -\dfrac{d\alpha(s)}{ds}, & s \in (0,\tau],\\
  \int_{-\tau}^0d\zeta^T(t,0)\alpha(-t),& s=0.
  \end{cases}
\end{equation*}
%where $\zeta^T$ denotes the transpose of $\zeta.$
\newline Note that, the domains of $\mathcal{A}$ and $\mathcal{A}^*$ are $C^1[-\tau,0]$ and $C^1[0,\tau]$ respectively. As
% \begin{equation*}
  $\mathcal{A}q(\theta) = \lambda(0)q(\theta)$
% \end{equation*}
$\bar{\lambda}(0)$ is an eigenvalue for $\mathcal{A}^*$, and
% \begin{equation*}
  $A^*q^* = -i\omega_0q^*$
% \end{equation*}
for some nonzero vector $q^*$. For $\phi \in C[-\tau,0]$ and $\psi \in C[0,\tau]$, define a bilinear inner product
\begin{equation}\label{inner_pro}
  \varsigma \langle \psi,\phi\rangle = \bar{\psi}(0).\phi(0)-\int_{\theta=-\tau}^0\int_{\varsigma=0}^\theta\bar\psi^T(\varsigma-\theta)d\zeta(\theta)\phi(\varsigma)d\varsigma.
\end{equation}
Then, $\langle\psi,A \phi\rangle = \langle A^*\psi,\phi\rangle $ for $\phi \in$ Dom$(\mathcal{A}),\psi \in$ Dom$(\mathcal{A}^*)$. Let $q(\theta) = e^{i\omega_0\theta}$ and $q^*(s) = De^{i\omega_0 s}$ be the eigenvectors for $ \mathcal{A}$ and $\mathcal{A}^*$ corresponding to the eigenvalues $+i\omega_0$ and $-i\omega_0$. Value of $D$ can be evaluated using (\ref{inner_pro}) and the relation $\langle q^*,q\rangle = 1$ as
\[
\begin{aligned}
 \langle q^*,q\rangle =&~  \bar{D}-\bar{D}\eta\int_{\theta = -\tau}^0\theta e^{i\omega_0\theta}\big(\xi_y\delta(\theta+\tau)\big)\mathrm{d}\theta \\
 \Rightarrow 1 =&~  \bar{D}+\bar{D}\eta\left(\tau\xi_ye^{-i\omega_0\tau}\right) \\
 \Rightarrow D = &~ \frac{1}{1+\eta\tau\xi_ye^{i\omega_0\tau}}. 
\end{aligned}
\]
%
% Again, using (\ref{inner_pro}) we show that $\langle q^*,\bar{q}\rangle=0$ as
% \[
% \begin{aligned}\normalsize
%   \langle q^*,\bar{q}\rangle = &~ \bar{D}+\frac{\bar{D}\eta}{2i\omega_0}\int_{\theta = -\tau}^0(e^{-i\omega_0\theta}-e^{i\omega_0\theta})\nonumber\\
%   %
%   &~~ \times \big(\xi_y\delta(\theta+\tau)+\xi_z\delta(\theta+\tau_2)\big)\mathrm{d}\theta,\nonumber\\
%   %
%   =&~ \bar{D}+\frac{\bar{D}\eta}{2i\omega_0}\big(\xi_y(e^{i\omega_0\tau}-e^{-i\omega_0\tau})\nonumber\\
%   %
%   &~ + \ \xi_z(e^{i\omega_0\tau_2}-e^{-i\omega_0\tau_2})\big),\nonumber\\
%   %
%   = &~ 0.\nonumber
% \end{aligned}
% \] 
Now we define
\begin{equation}
\begin{aligned}
z(t) &= \langle q^*.u_t\rangle,  \\
w(t,\theta)&= u_t(\theta)-2Re\{z(t)q(\theta)\}.
\end{aligned}
\end{equation}
Then,  on the center manifold $C_0$, $w(t,\theta)=w\big(z(t), \bar{z}(t),\theta\big)$ where
\begin{equation}\label{w_equ}
  w(z,\bar{z},\theta)=w_{20}(\theta)\frac{z^2}{2}+w_{11}(\theta)z\bar{z}+w_{02}(\theta)\frac{\bar{z}^2}{2}+\cdots.
\end{equation}
In effect, $z$ and $\bar{z}$ are local coordinates for manifold in $C$ in the directions
of $q^*$ and $\bar{q}^*$, respectively. The existence of the center
manifold $C_0$ enables us to reduce (\ref{eq:auto_noneq_matrix}) to an ordinary differential
equation for a single complex variable on $C_0$. At $\mu=0$, we have
\begin{eqnarray}
z'(t) &\hspace{-2mm}=&\hspace{-2mm} \langle q^*,\mathcal{A}y_t+\mathcal{R}u_t\rangle\nonumber\\
&\hspace{-2mm}=&\hspace{-2mm} i\omega_0z(t)+\bar{q}^*(0).\mathcal{F}\big(w(z,\bar{z},\theta)+2\mathbf{Re}\{z(t)q(\theta)\}\big)\nonumber\\
&\hspace{-2mm}=&\hspace{-2mm} i\omega_0z(t)+\bar{q}^*(0).\mathcal{F}_0(z,\bar{z})\label{z_dot_prod}
\end{eqnarray}
which can be written as
\begin{equation}\label{abb_z}
  z'(t) = i\omega_0z(t)+g(z,\bar{z}).
\end{equation}
% The next objective is to expand $g$ in powers of $z$  and $\bar{z}$ .
% However, we also have to determine the coefficients $w_{ij}(\theta)$ in (\ref{w_equ}).
% Once the $w_{ij}$ have been determined, the differential equation (\ref{z_dot_prod}) for $z$ would be explicit [as abbreviated in (\ref{abb_z})] where 
Expanding $g(z,\bar{z})$ in powers of $z$ and $\bar{z}$, we get
\begin{eqnarray}
  g(z,\bar{z}) &=& \bar{q}^*(0).\mathcal{F}_0(z,\bar{z})\nonumber\\
  &=& g_{20}\frac{z^2}{2}+g_{11}z\bar{z}+g_{02}\frac{\bar{z}^2}{2}+g_{21}\frac{z^2\bar{z}}{2}+\cdots.\nonumber
\end{eqnarray}
Following \cite{hassard1981}, we write
\begin{equation}\label{wdash}
  w' = u_t'-z'q-\bar{z}'\bar{q}.
\end{equation}
From (\ref{eq:auto_noneq_matrix}) and (\ref{abb_z}), we get
\begin{equation*}
  w' =
    \begin{cases}
  Aw-2\mathbf{Re}\{\bar{q}^*(0).\mathcal{F}_0q(\theta)\}, & \theta \in [-\tau_2,0)\\
  Aw-2\mathbf{Re}\{\bar{q}^*(0).\mathcal{F}_0q(0)\}+\mathcal{F}_0,& \theta=0
  \end{cases}
\end{equation*}
which can be written as
\begin{equation}\label{eq:matrix_w}
  w' = Aw+H(z,\bar{z},\theta),
\end{equation}
using (\ref{abb_z}), where
\begin{equation}\label{H_equ}
  H(z,\bar{z},\theta)=H_{20}(\theta)\frac{z^2}{2}+H_{11}(\theta)z\bar{z}+H_{02}(\theta)\frac{\bar{z}^2}{2}+\cdots.
\end{equation}
Now, on $C_0$, near the origin
% \begin{equation*}
  $w' = w_zz'+w_{\bar{z}}\bar{z}'$.
%  \end{equation*}
Using (\ref{w_equ}) and (\ref{abb_z}) to replace $w_z,z'$, and equating this with (\ref{eq:matrix_w}), we get
\begin{eqnarray}
  (2i\omega_0-A)w_{20}(\theta) &=& H_{20}(\theta)\label{eq1} \\
  -Aw_{11}(\theta) &=& H_{11}(\theta)\label{eq2}\\
  (2i\omega_0-A)w_{02}(\theta) &=& H_{02}(\theta)\label{eq2.1}.
\end{eqnarray}
as in \cite{hassard1981}. From (\ref{wdash}), we get
\begin{small}
\begin{align}
  u_t(\theta) &= w(z,\bar{z},\theta)+zq(\theta)+\bar{z}\bar{q}(\theta) \nonumber\\
   &= w_{20}(\theta)\frac{z^2}{2}+w_{11}z\bar{z}+w_{02}(\theta)\frac{\bar{z}^2}{2}+ze^{i\omega_0\theta}+\bar{z}e^{-i\omega_0\theta}+\cdot\cdot\nonumber
\end{align}
\end{small}
from which $u_t(0)$ and $u_t(-\tau)$ can be determined. 
As we only require the coefficients of $z^2,z\bar{z},\bar{z}^2$ and $z^2\bar{z}$, we have

\begin{eqnarray}
%   u_t^2(0)  &\hspace{-1.5mm}=&\hspace{-1.5mm} \big(w(z,\bar{z},0)+z+\bar{z}\big)^2 \nonumber\\
%    &\hspace{-1.5mm}=&\hspace{-1.5mm} z^2+\bar{z}^2+2z\bar{z}+z^2\bar{z}\big(2w_{11}(0)\nonumber\\
%    &\hspace{-1.5mm}&\hspace{-1.5mm}+ w_{20}(0)\big)+\cdots.\nonumber\\
u_t(0)u_t(-\tau)&\hspace{-1.5mm}=&\hspace{-1.5mm} \big(w(z,\bar{z},0)+z+\bar{z}\big)\nonumber\\
&\hspace{-1.5mm}&\hspace{-1.5mm} \times \big(w(z,\bar{z},-\tau)+ze^{-i\omega_0\tau}+\bar{z}e^{i\omega_0\tau}\big)   \nonumber\\
&\hspace{-1.5mm}=&\hspace{-1.5mm} z^2e^{-i\omega_0\tau}+z\bar{z}(e^{i\omega_0\tau}+e^{-i\omega_0\tau})+ \bar{z}^2e^{i\omega_0\tau}\nonumber\\
&\hspace{-1.5mm}&\hspace{-1.5mm} +\ z^2\bar{z}\left(w_{11}(0)e^{-i\omega_0\tau}+\frac{w_{20}(0)}{2}e^{i\omega_0\tau}\right.\nonumber\\
&\hspace{-1.5mm}&\hspace{-1.5mm}\left.+\ w_{11}(-\tau)+\frac{w_{20}(-\tau)}{2}\right)+\cdots\ .\nonumber\\
%&\hspace{-1.5mm}&\hspace{-1.5mm}+\cdots, \forall i = 1,2.\nonumber\\
%**********************
%***********
% \end{eqnarray}
% \begin{eqnarray}
u_t^2(-\tau) &\hspace{-1.5mm}=&\hspace{-1.5mm} \big(w(z,\bar{z},-\tau)+ze^{-i\omega_0\tau}+\bar{z}e^{i\omega_0\tau}\big)^2\nonumber\\
&\hspace{-1.5mm}=&\hspace{-1.5mm} z^2e^{-2i\omega_0\tau}+\bar{z}^2e^{2i\omega_0\tau}+2z\bar{z}\nonumber\\   
&\hspace{-1.5mm}&\hspace{-1.5mm}+\ z^2\bar{z}\left(2e^{-i\omega_0\tau}w_{11}(-\tau)\right.\nonumber\\
&\hspace{-1.5mm}&\hspace{-1.5mm}\left.+\ e^{i\omega_0\tau} w_{20}(-\tau)\right) +\cdots \ .\nonumber\\
%*******************
%************************
% u_t^3(0)&\hspace{-2mm}=&\hspace{-2mm}\big(w(z,\bar{z},0)+z+\bar{z}\big)^3\nonumber\\
% &\hspace{-2mm}=&\hspace{-2mm} 3z^2\bar{z}+\cdots.\nonumber\\
%************************
u_t^2(0)u(-\tau) &\hspace{-2mm}=&\hspace{-2mm}  (w(z,\bar{z},0)+z+\bar{z})^2\nonumber\\
&\hspace{-2mm}&\hspace{-2mm}\times (w(z,\bar{z},-\tau)+ze^{-i\omega_0\tau}+\bar{z}e^{i\omega_0\tau})\nonumber\\
&\hspace{-2mm}=&\hspace{-2mm} z^2\bar{z}(2e^{-2i\omega_0\tau}+e^{i\omega_0\tau})+\cdots.\nonumber\\
%************************
u_t(0)u^2(-\tau) &\hspace{-2mm}=&\hspace{-2mm}  \big(w(z,\bar{z},0)+z+\bar{z}\big)\nonumber\\
&\hspace{-2mm}&\hspace{-2mm}\times \big(w(z,\bar{z},-\tau)+ze^{-i\omega_0\tau}+\bar{z}e^{i\omega_0\tau}\big)^2\nonumber\\
&\hspace{-2mm}=&\hspace{-2mm} z^2\bar{z}(e^{-2i\omega_0\tau}+2)+\cdots. \nonumber\\
%**********************
u_t^3(-\tau) &\hspace{-2mm}=&\hspace{-2mm}  \big(w(z,\bar{z},-\tau)+ze^{-i\omega_0\tau}+\bar{z}e^{i\omega_0\tau}\big)^3\nonumber\\
&\hspace{-2mm}=&\hspace{-2mm} 3z^2\bar{z}e^{-i\omega_0\tau}+\cdots. \nonumber
\end{eqnarray}

%%%%%%%%%%%%%%%%%%%%%%%%%%%%%%%%%%%%%%%%%%%%%%%%%%%%%%%%%%%%%%%%%%%%%%%%%%%%%%%%%%%%%%%%%%%%%%%%%%%%%%%%%%%%%%%%%%%%%%%%%%%%%%%%%%
Recall that
\begin{eqnarray}
   g(z,\bar{z}) &=& \bar{q}^*(0).\mathcal{F}_0(z,\bar{z})\nonumber\\
   &=& g_{20}\frac{z^2}{2}+g_{11}z\bar{z}+g_{02}\frac{\bar{z}^2}{2}+g_{21}\frac{z^2\bar{z}}{2}+\cdots.\nonumber
\end{eqnarray}
Comparing the coefficients of $z^2,z\bar{z},\bar{z}^2$, and $z^2\bar{z}$, we get
\begin{eqnarray}
g_{20}=&&\hspace{-6mm}\bar{D}\eta[2\xi_{xx}+2\xi_{xy}e^{-i\omega_0\tau}+2\xi_{yy}e^{-2i\omega_0\tau}]\nonumber\\
g_{11}=&&\hspace{-6mm}\bar{D}\eta[2\xi_{xx}+\xi_{xy}(e^{-i\omega_0\tau}+e^{i\omega_0\tau})+2\xi_{yy}]\nonumber\\
g_{02}=&&\hspace{-6mm}\bar{D}\eta[2\xi_{xx}+2\xi_{xy}e^{i\omega_0\tau}+2\xi_{yy}e^{2i\omega_0\tau}]\nonumber\\
g_{21}=&&\hspace{-6mm}\bar{D}\eta[2\xi_{xx}\big(2w_{11}(0)+w_{20}(0)\big)+\xi_{xy}\big(2w_{11}(0)e^{-i\omega_0\tau}\nonumber\\
%&&\hspace{-6mm}+\ \xi_{xy}\big(2w_{11}(0)e^{-i\omega_0\tau}+ w_{20}(0)e^{i\omega_0\tau}\nonumber\\
&&\hspace{-6mm}+ w_{20}(0)e^{i\omega_0\tau}+2w_{11}(-\tau)+w_{20}(-\tau)\big)\nonumber\\
&&\hspace{-6mm}+\ \xi_{xz}\big(2w_{11}(0)e^{-i\omega_0\tau_2}+w_{20}(0)e^{i\omega_0\tau_2}\nonumber\\
&&\hspace{-6mm}+\ 2w_{11}(-\tau_2)+w_{20}(-\tau_2)\big)\nonumber\\
&&\hspace{-6mm}+\ \xi_{yy}\big(4w_{11}(-\tau)e^{-i\omega_0\tau}+2w_{20}(-\tau)e^{i\omega_0\tau}\big)\nonumber\\
&&\hspace{-6mm}+\ 6\xi_{xxx}+\xi_{xyy}(2e^{-2i\omega_0\tau}+4) \nonumber\\
&&\hspace{-6mm}+\ \xi_{xxy}(2e^{i\omega_0\tau}+4e^{-i\omega_0\tau})+6\xi_{yyy}e^{-i\omega_0\tau}].\label{g21}
\end{eqnarray}
% Substituting the values, we get\\
% \begin{eqnarray}
% g_{20}=&&\hspace{-6mm}-4j\bar{D}\eta\xi_{xy}sin{(\omega_0\tau)}.\label{g20}\\
% g_{11}=&&\hspace{-6mm}0.\label{g11}\\
% g_{02}=&&\hspace{-6mm}4j\bar{D}\eta\xi_{xy}sin{(\omega_0\tau)}..\label{g02}
% \end{eqnarray}
For $\theta \in [-\tau,0)$, we have
\begin{small}
\begin{eqnarray*}
H(z,\bar{z},\theta) &&\hspace{-6mm}= -2Re\{\bar{q}^*(0).\mathcal{F}_0q(\theta)\}=-g(z,\bar{z})q(\theta) -\bar{g}(z,\bar{z})\bar{q}(\theta)\\
&&\hspace{-6mm}=-\Big(g_{20}\frac{z^2}{2}+g_{11}z\bar{z} + g_{02}\frac{\bar{z}^2}{2}+\cdots\Big)q(\theta)\\
&&\hspace{-6mm} \ \ \ -\Big(\bar{g}_{20}\frac{\bar{z}^2}{2}+\bar{g}_{11}z\bar{z} + \bar{g}_{02}\frac{z^2}{2}+\cdots\Big)\bar{q}(\theta).
\end{eqnarray*}
\end{small}
Now using (\ref{H_equ}), we obtain
\begin{equation*}
  H_{20}(\theta) = -g_{20}q(\theta)-\bar{g}_{02}\bar{q}\theta,\ \quad  H_{11}(\theta) = -g_{11}q(\theta)-\bar{g}_{11}\bar{q}\theta.
\end{equation*}
From (\ref{eq:Atheta}), (\ref{eq1}) and (\ref{eq2}), we derive
\begin{eqnarray*}
  w'_{20}(\theta) =&&\hspace{-6mm} 2i\omega_0w_{20}(\theta)+g_{20}q(\theta)+\bar{g}_{02}\bar{q}(\theta), \label{eq4}\\
  w'_{11}(\theta) =&&\hspace{-6mm} g_{11}q(\theta)+\bar{g}_{11}\bar{q}(\theta).\label{eq5}
\end{eqnarray*}
Solving the above differential equations yields
\begin{eqnarray}
w_{20}(\theta) &\hspace{-3mm}=&\hspace{-3mm} -\frac{g_{20}}{i\omega_0}q(0)e^{i\omega_0\theta} -\frac{\bar{g}_{02}}{3i\omega_0}\bar{q}(0)e^{-i\omega_0\theta}+Ee^{2i\omega_0\theta}\nonumber\\&&\label{w20}\\
w_{11}(\theta) &\hspace{-3mm}=&\hspace{-3mm} \frac{g_{11}}{i\omega_0}q(0)e^{i\omega_0\theta} -\frac{\bar{g}_{11}}{i\omega_0}\bar{q}(0)e^{-i\omega_0\theta}+F\label{w11}
\end{eqnarray}
for some $E$ and $F$. For $\theta=0$, we get
\begin{eqnarray}
  H(z,\bar{z},0)=&&\hspace{-6mm} -2Re(\bar{q}^*.\mathcal{F}_0q(0))+\mathcal{F}_0,\nonumber\\
  H_{20}(0) =&&\hspace{-6mm} -g_{20}q(0) -\bar{g}_{02}\bar{q}(0)\nonumber\\
  &&\hspace{-6mm}+\eta\big[2\xi_{xx}+2\xi_{xy}e^{-i\omega_0\tau}+2\xi_{yy}e^{-2i\omega_0\tau}\big]\nonumber \label{H_20}\\
  H_{11}(0) =&&\hspace{-6mm} -g_{11}q(0) -\bar{g}_{11}\bar{q}(0) \nonumber\\
  &&\hspace{-6mm}+\eta\big[2\xi_{xx}+\xi_{xy}(e^{-i\omega_0\tau}+e^{i\omega_0\tau})+ 2\xi_{yy}\big].\nonumber \label{H_11}
\end{eqnarray}
Using (\ref{eq:Atheta}), (\ref{eq1}) and (\ref{eq2}), we get
\begin{eqnarray}
  &&\eta\xi_yw_{20}(-\tau)+\eta\xi_xw_{20}(0)-2i\omega_0w_{20}(0)\nonumber\\
  &&\quad=g_{20}q(0)+\bar{g}_{02}\bar{q}(0)\nonumber\\
  &&\quad\quad -\eta\big[2\xi_{xx}+2\xi_{xy}e^{-i\omega_0\tau}+\ 2\xi_{yy}e^{-2i\omega_0\tau}\big]\label{eq8} \\
  &&\eta\xi_yw_{11}(-\tau)+\eta\xi_xw_{11}(0)\nonumber\\
  &&\quad=  g_{11}q(0)+\bar{g}_{11}\bar{q}(0)\nonumber\\
  &&\quad\quad  -\eta\big[2\xi_{xx}+\xi_{xy}(e^{-i\omega_0\tau}+e^{i\omega_0\tau})+ 2\xi_{yy}\big]. \label{eq9}
\end{eqnarray}
% We have the solution for $w_{20}(\theta)$ and $w_{11}(\theta)$ from (\ref{w20}) and (\ref{w11}) respectively. 
% Hence, evaluate $w_{11}(0)$, $w_{20}(0)$, $w_{11}(-\tau)$, $w_{20}(-\tau)$,  $w_{11}(-\tau_2)$, $w_{20}(-\tau_2)$, substitute into (\ref{eq8}) and (\ref{eq9}) respectively, and calculate $E,F$ as
Evaluate $w_{11}(0)$, $w_{20}(0)$, $w_{11}(-\tau)$ and $w_{20}(-\tau)$ using (\ref{w20}) and (\ref{w11}), and substituting in (\ref{eq8}) and (\ref{eq9}), we get $E$ and $F$ as\\
\begin{small}
\begin{equation*}
   E = \dfrac{-g_{20}}{\bar{D}(\eta\xi_x+\eta\xi_ye^{-2i\omega_0\tau}-2i\omega_0)},\ \quad   F = \dfrac{-g_{11}}{\bar{D}\eta(\xi_y+\xi_x)}.
\end{equation*}
\end{small}
% where
% \begin{eqnarray}
% % \nonumber to remove numbering (before each equation)
%   \Phi_1 =&&\hspace{-6mm}(-2i\omega_0)\left(\frac{g_{20}}{i\omega_0}+\frac{\bar{g}_{02}}{3i\omega_0}\right) +\ \eta\xi_y\left(\frac{g_{20}}{i\omega_0}e^{-i\omega_0\tau}\right.\nonumber\\
%   &&\hspace{-6mm}\left.+\ \frac{\bar{g}_{02}}{3i\omega_0}e^{i\omega_0\tau}\right) +\eta\xi_z\left(\frac{g_{20}}{i\omega_0}e^{-i\omega_0\tau_2}\right.\nonumber\\
%   &&\hspace{-6mm}\left.+\frac{\bar{g}_{02}}{3i\omega_0}e^{i\omega_0\tau_2}\right)-H_{20}(0)\nonumber\\
%   \Phi_2 =&& \hspace{-6mm}-\ \eta\xi_y\left(\frac{g_{11}}{i\omega_0}e^{-i\omega_0\tau}\right.\nonumber\\
%   &&\hspace{-6mm}\left.-\ \frac{\bar{g}_{11}}{i\omega_0}e^{i\omega_0\tau}\right)- \eta\xi_z\left(\frac{g_{11}}{i\omega_0}e^{-i\omega_0\tau_2}\right.\nonumber\\
%   &&\hspace{-6mm}\left.-\ \frac{\bar{g}_{11}}{i\omega_0}e^{i\omega_0\tau_2}\right)-H_{11}(0)\nonumber.
% \end{eqnarray}

Thus, the stability of the bifurcating solutions can now be investigated using \cite{hassard1981}. The quantities required to study the nature of the Hopf bifurcation are as follows\\
\begin{eqnarray}
\mu_2 &&\hspace{-6mm}= \dfrac{-\operatorname{Re}[c_1(0)]}{\alpha'(0)},\quad \beta_2= 2\operatorname{Re}[c_1(0)],\nonumber
%\beta &&\hspace{-6mm}= \epsilon^2\beta_2+\mathcal{O}(\epsilon^4)\quad \beta_2 = 2\operatorname{Re}[c_1(0)]\quad \epsilon = \sqrt{\frac{\mu}{\mu_2}}\nonumber\\\label{49}
\end{eqnarray}
where $\ \alpha'(0) =\mathbf{Re}({d\lambda}/{d\eta})_{\eta=\eta_c},$
%\begin{small}
\begin{equation*}
c_1(0) = \dfrac{i}{2\omega_0}\left(g_{20}g_{11}-2|g_{11}|^2-\dfrac{1}{3}|g_{02}|^2\right)+\dfrac{g_{21}}{2}.
\end{equation*}
%\end{small}
% $\alpha'(0)$ is the real component of $(d\lambda/d\eta)$ evaluated at $\eta = \eta_c$
The direction and stability of the Hopf bifurcation is determined by the sign of $\mu_2$ and  $\beta_2$, respectively.
If $\mu_2 >0(\mu_2 <0)$ then the Hopf bifurcation is \emph{super-critical(sub-critical)}. Similarly, the bifurcating solutions are \emph{asymptotically orbitally stable(unstable)} if $\beta_2<0(\beta_2>0)$.\\

Using the calculations outlined above, we obtain the expression for $\mu_2$ as 
\begin{scriptsize}
 \begin{eqnarray}
 \mu_2 =&& \hspace{-5mm} \dfrac{1}{b^2 (1+\epsilon) (1-\epsilon^2) \hat{\epsilon} (5-4\epsilon)} \Bigg\{
%  \nonumber \\
%  && \hspace{-5mm} 
 \xi^2_{xx} \bigg(\check{\epsilon}(12\epsilon+8) + \hat{\epsilon} (8\epsilon^2 -18\epsilon +4)\bigg)\nonumber \\
 && \hspace{-8mm} +\xi^2_{xy} \bigg(\check{\epsilon}(4\epsilon^3 - 14\epsilon^2 + 11\epsilon - 1) + \hat{\epsilon} (-8\epsilon^3 + 12\epsilon^2 - 7\epsilon +3)\bigg) \nonumber \\
  && \hspace{-8mm} +\xi^2_{yy} \bigg(\check{\epsilon}(-8\epsilon^3 -8\epsilon^2 + 26\epsilon - 4) + \hat{\epsilon} (-4\epsilon^2 -12\epsilon +22)\bigg) \nonumber \\
   && \hspace{-8mm} +\xi_{xy} \xi_{xx} \bigg(\check{\epsilon}(-18\epsilon^2 + 33\epsilon -9) + \hat{\epsilon} (-8\epsilon^3 + 26\epsilon^2 -19\epsilon +7)\bigg) \nonumber \\
    && \hspace{-8mm} +\xi_{xy} \xi_{yy} \bigg(\check{\epsilon}(8\epsilon^4 + 8\epsilon^3 - 32\epsilon^2 + 19\epsilon - 9) + \hat{\epsilon} (4\epsilon^3 + 20\epsilon^2 - 37\epsilon + 7)\bigg)\nonumber \\
     && \hspace{-8mm} +\xi_{xx} \xi_{yy} \bigg(\check{\epsilon}(-12\epsilon^2 +30\epsilon - 18) + \hat{\epsilon} (16\epsilon^2 -30\epsilon +14)\bigg) \Bigg\} \nonumber \\
     && \hspace{-8mm} + \dfrac{1}{b (1-\epsilon^2) \hat{\epsilon}} \Bigg\{ \xi_{xxx} \bigg(-3\check{\epsilon} - \hat{\epsilon} \epsilon \bigg) +  \xi_{xyy} \bigg(-\check{\epsilon} (1+2\epsilon^2)-3\hat{\epsilon} \epsilon \bigg)\nonumber \\   && \hspace{-8mm} + \xi_{xxy} \bigg(3\check{\epsilon} \epsilon + \hat{\epsilon} (1+2\epsilon^2) \bigg) +  \xi_{yyy} \bigg(3\check{\epsilon} \epsilon + 3\hat{\epsilon} \bigg)\Bigg\} \label{eq:mu2_gen_expression}
\end{eqnarray}
\end{scriptsize}
where $\epsilon = a/b = \xi_x/\xi_y$, $\check{\epsilon} = \sqrt{1 - \epsilon^2}$, and $\hat{\epsilon} = \cos^{-1}(-\epsilon)$. It is to be noted that the result that we obtained is not just confined to the chemical oscillator but can also be extended to various first-order non-linear delay dynamical systems.

% The results of our study can also be used to analyze the bifurcation properties of the Nicholson's Blowflies equation, which has been extensively used in the context of population dynamics. 
\textbf{Example 1}. Consider the following equation:
\begin{equation}
 \dot{N}(t) = -\gamma N(t) + p N(t - \tau) e^{- N(t-\tau)/x_0}.\label{eq:nicholsonbfeqn}
\end{equation}
The above equation is called the \textit{Nicholson's blowflies equation} \cite{nicholsonbfref}, which has a variety of applications in the context of population dynamics. Here, $N(t)$ is the size of the population at time $t$, $p$ is the maximum per capita daily egg production rate, $x_0$ is the size at which the population reproduces at the maximum rate, $\gamma$ is the per capita daily adult death rate, and $\tau$ is the generation time. In the literature, there are many studies on the stability and oscillations of \eqref{eq:nicholsonbfeqn}, for example, see \cite{feng2002_nic,gyori2002_nic,saker2002_nic}. However, there exists little research on the Hopf bifurcation properties of \eqref{eq:nicholsonbfeqn}. In \cite{shi2002_nic}, the  impact of loss of stability of \eqref{eq:nicholsonbfeqn} has been analyzed, but, it is for a particular choice of parameter values. To rule out being mislead by a particular choice, 
a detailed Hopf bifurcation analysis is required. To that end, we now analyze the Hopf bifurcation of \eqref{eq:nicholsonbfeqn} using the general result \eqref{eq:mu2_gen_expression} obtained in this study.

% To our knowledge, the bifurcation properties of \eqref{eq:nicholsonbfeqn} have not been intensively studied in the literature. 
Including the exogenous bifurcation parameter ($\eta$) in (\ref{eq:nicholsonbfeqn}), and expanding it using Taylor series, we obtain
\begin{equation}
 \dot{u}(t) =\eta \big(-a u(t) - b u(t-\tau) + \xi_{yy}u^2(t - \tau)  %\notag \\
%  &+ \xi_{xy}u(t)u(t-\tau)+ \xi_{yy}u^2(t-\tau) + \xi_{xxx}u^3(t) \notag  \\
% 		  & + \xi_{xxy}u^2(t)u(t-\tau)  + \xi_{xyy}u(t)u^2(t-\tau) \notag  \\
		   + \xi_{yyy}u^3(t-\tau) \big),
\label{eq:taylor_nicholson}
		  \end{equation}
where $a = -\xi_x = \gamma$,\  $b = -\xi_y = -\gamma (1 - \ln(p/\gamma) )$, \  $\xi_{yy} = -\frac{\gamma}{x_0}(2 - \ln(p/\gamma))$, and  $\xi_{yyy} = \frac{\gamma}{x^2_0}(3 - \ln(p/\gamma))$.\\

The Hopf condition is: $\eta_c \tau = {\cos^{-1}(-a/b)}/{\sqrt{b^2 - a^2}}$, where $\eta_c$ is the critical value of $\eta$ which induces a Hopf bifurcation. Using \eqref{eq:mu2_gen_expression}, we obtain $\mu_2$ for \eqref{eq:taylor_nicholson} as
 \begin{scriptsize}
 \begin{align}
 \mu_2 = & \ \dfrac{\xi^2_{yy} \Big({\sqrt{1 - \epsilon^2}(-8\epsilon^3 -8\epsilon^2 + 26\epsilon - 4) - \cos^{-1}(-\epsilon) (4\epsilon^2 +12\epsilon -22)} \Big)}{b^2 (1+\epsilon) (1-\epsilon^2) \cos^{-1}(-\epsilon) (5-4\epsilon)}  \nonumber \\
  & +  \dfrac{\xi_{yyy}\Big(3\epsilon \sqrt{1 - \epsilon^2}  + 3\cos^{-1}(-\epsilon)\Big)}{{b (1-\epsilon^2) \cos^{-1}(-\epsilon)}}, \label{eq:mu2_nicholson1}
\end{align}
\end{scriptsize}
where $\epsilon =  \frac{\xi_x}{\xi_y} = \frac{1}{\ln(P/\gamma) -1}$. On simplification, we get
\begin{scriptsize}
\begin{align}
\label{eq:nicholsonmu2_final}
\mu_2 =& \frac{1}{x^2_0} \Bigg\{ \dfrac{(1-\epsilon)}{ (1+\epsilon)^2 \cos^{-1}(-\epsilon) (5 - 4\epsilon)}  \bigg( \sqrt{1 - \epsilon^2}(-8\epsilon^3 -8\epsilon^2 + 26\epsilon - 4) \\ \nonumber
& + \cos^{-1}(-\epsilon) (-4\epsilon^2 -12\epsilon +22) \bigg) + \frac{2\epsilon - 1}{(1 - \epsilon^2) \cos^{-1}(-\epsilon)} \bigg(3\epsilon\sqrt{1-\epsilon^2}\\ \nonumber
& + 3\cos^{-1}(-\epsilon) \bigg) \Bigg\}. 
\end{align}
\end{scriptsize}
\begin{figure}[hbtp!]   
\vspace{-10mm}
\centering
%\mypsfrag{-0.6}{-0.2}\mypsfrag{0.0}{0.2}
\psfrag{0}{  \hspace{-2mm} $0$}
\psfrag{50}{ \small \hspace{-2mm} $50$}
\psfrag{100}{ \small  \hspace{-3mm}  $100$}
\psfrag{150}{ \small \hspace{-2mm}  $150$}
\psfrag{0.0}{ \small \hspace{-0mm} $0$}
\psfrag{0.2}{ \small \hspace{-1mm} $0.2$}
\psfrag{1.0}{ \small \hspace{-1mm} $1.0$}
\psfrag{0.4}{\small \hspace{-1mm} $0.4$}
\psfrag{0.6}{ \small \hspace{-1mm} $0.6$}
\psfrag{0.8}{ \small \hspace{-1mm} $0.8$}
\psfrag{pi/4}{ \small $\pi/4$}
\psfrag{pi/2}{ \small $\pi/2$}
\psfrag{x}{  $\epsilon$}
\psfrag{myu2}{  $\mu_2$}
\psfrag{0.8}{ \small $0.8$}
%\psfrag{M}[b][t]{\hspace{0mm}\footnotesize $\mu_2$}
\psfrag{M}{\hspace{0mm}\footnotesize $\mu_{2}$}
\psfrag{B}{\hspace{0mm}\footnotesize $\beta_{2}$}
\includegraphics[height=2.15in,width=2.97in]{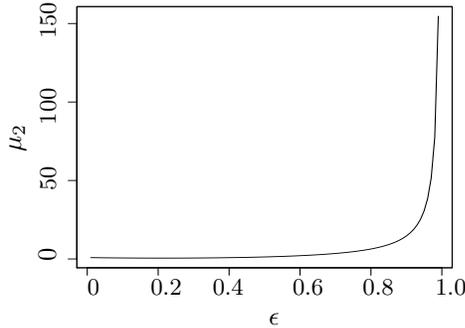}
\vspace{-3mm}
\caption{ The plot of $\mu_2$ Vs $\epsilon$. Observe that, for all values of $\epsilon \in (0,1)$, the value of $\mu_2$ is positive, which implies that the Hopf bifurcation is super-critical.}\label{fig:nicholsonmyu2plot}
\end{figure}
From \eqref{eq:nicholsonmu2_final}, we can see that the value of $x_0$ does not affect the sign of $\mu_2$, and hence we consider $x_0 =1$, and plot $\mu_2$ for $\epsilon \in (0,1)$. From Figure \ref{fig:nicholsonmyu2plot}, we can observe that the value of $\mu_2$ is always greater than zero. This implies that the type of Hopf bifurcation is super-critical, which would lead to stable and small-amplitude limit cycles.

\balance

\end{document}